\documentclass[a4paper,10pt,oneside]{article}      
\usepackage[english]{babel}
\usepackage[utf8]{inputenc}
\usepackage[T1]{fontenc}
\usepackage{amssymb,amsthm,amsmath}
\usepackage{verbatim}
\usepackage{enumitem} 
\usepackage{a4wide}
\usepackage{hyperref}
\usepackage{definitions}
\usepackage[sorting=nyt, style=alphabetic, firstinits=true, doi=false,isbn=false,url=false, backend=biber]{biblatex}
\bibliography{density-waring}
\AtEveryBibitem{\ifentrytype{book}{\clearfield{pages}}}

\usepackage{xcolor}
\hypersetup{
    colorlinks,
    linkcolor={red!50!black},
    citecolor={blue!50!black},
    urlcolor={blue!80!black}
}

\makeatletter
\newcommand{\mynameis}[1]{#1\renewcommand{\@currentlabel}{#1}}
\makeatother

\allowdisplaybreaks

\begin{document}

\title{A Density version of Waring's problem}
\author{Juho Salmensuu}
\date{}
\maketitle

\begin{abstract}
   In this paper, we study a density version of Waring's problem. We prove that a positive density subset of $k$th-powers forms an asymptotic additive basis of order $O(k^2)$ provided that the relative lower density of the set is greater than $(1 - \mathcal{Z}_k^{-1}/2)^{1/k}$, where $\mathcal{Z}_k$ is a certain constant depending on $k$ for which it holds that $\mathcal{Z}_k > 1$ for every $k$ and $\lim_{k \rightarrow \infty} \mathcal{Z}_k = 1$. 
\end{abstract}

\section{Introduction}

\subsection{Statements of results}\label{section_statement_of_results}

In this paper, we investigate, when a positive density subset of $k$th-powers forms an asymptotic additive basis. This problem is motivated by similar results related to Goldbach's problem \cite{li-pan}, \cite{shao}. For example Shao \cite{shao} proved that if $A$ is a subset of the primes, and the lower density of $A$ in the primes is larger than $5/8$, then all sufficiently large odd positive integers can be written as the sum of three primes in $A$. The key to studying these kinds of problems is the transference principle introduced by Green \cite{green-transference}.

Let $k \geq 2$ be an integer. Set $\N^{(k)} := \{t^k : t \in \N\}$ and $\Z_m^{(k)} := \{t^k : t \in \Z_m\}$, where $\Z_m := \Z/m\Z$. Let $A \subseteq \N^{(k)} $. Define
\begin{displaymath}
\delta_A = \underline{\delta}(A) := \liminf_{N \rightarrow \infty} \frac{|A \cap [N]|}{|\N^{(k)} \cap [N]|},
\end{displaymath}
where $[N] := \{1, \dots, N\}$. For $m \in \N$ let $P(m) := \prod_{p \leq m} p^k$ and
\begin{equation} \label{defn_Zk}
\mathcal{Z}_k :=  \lim_{m \rightarrow \infty} \frac{|\Z_{P(m)}^{(k)}|}{|\{a \in \Z_{P(m)}^{(k)} \mid (a, P(m)) = 1\}|}.
\end{equation}
We prove later that $\lim_{k \rightarrow \infty} \mathcal{Z}_k = 1$.

For each prime $p$ and $k \in \N$, define $\tau(k, p)$ so that $p^{\tau(k, p)} || k$, where $p^{h} || k$ means that $p^{h} | k$ and $p^{h+1} \nmid k$. Let 
\begin{equation} \label{defn_R_k}
R_k := \prod_{(p-1) | k}p^{\eta(k, p)},
\end{equation}
where 
\begin{equation} \label{defn_eta}
\eta(k, p) := \left\{
	\begin{array}{ll}
		\tau(k, p) + 2 & \mbox{if } p=2 \text{ and } \tau(k, p) > 0\\
		\tau(k, p) + 1 & \mbox{otherwise}
	\end{array}
\right.
\end{equation}
For $n \in \N$ let the function $\omega(n)$ denote the number of distinct prime divisors of $n$.

Our main result is the following.

\begin{Theorem} \label{theorem_main} Let $s, k \in \N$, $k \geq 2$, $s > \max(16k\omega(k) + 4k + 3, k^2+k)$ and let $A \subseteq \N^{(k)}$ be such that $\underline{\delta}(A) > (1 - \mathcal{Z}_k^{-1}/2)^{1/k}$. Then, for all sufficiently large integers $n \equiv s \Mod{R_k}$, we have $n \in sA$.
\end{Theorem}

In the proof of the last theorem, due to some technical difficulties, we need to restrict to the elements of $A$ which do not have small prime factors. This leads to the congruence condition in the previous theorem. We expect that the congruence condition could be removed. In the following corollary, we have done so at the cost of some extra summands.

\begin{Corollary}  Let $s, k \in \N$, $k \geq 2$, $s > \max(16k\omega(k) + 4k + 3, k^2+k) + R_k$ and let $A \subseteq \N^{(k)}$ be such that $\underline{\delta}(A) > (1 - \mathcal{Z}_k^{-1}/2)^{1/k}$. Then, for all sufficiently large $n \in \N$, we have $n \in sA$.
\end{Corollary}

\begin{proof} For $q \in \N$, let $A_q := \{b \in \Z_{q} \mid \exists a \in A: b \equiv a \Mod{q}\}$. Let $P = \{p : p-1 | k\}$. By $\underline{\delta}(A)  > (1/2)^{1/k} > \frac{k}{k + 1} \geq \max_{p \in P} \frac{p-1}{p}$ and Fermat's little theorem we see that $A_{p^{\eta(k,p)}} = \{0, 1\}$ for all $p \in P$. Hence $sA_{p^{\eta(k,p)}} = \Z_{p^{\eta(k,p)}}$ for all $s \geq p^{\eta(k,p)}$ and $p \in P$. Therefore by the Chinese remainder theorem  $sA_{R_k} = \Z_{R_k}$ for all $s \geq \max_{p: p-1 | k} p^{\eta(k,p)}$. The rest now follows from Theorem \ref{theorem_main}.
\end{proof}

The density condition in Theorem \ref{theorem_main} is not optimal. We expect the result to hold as long as $A$ is not a subset of a non-trivial \footnote{By \textit{non-trivial} we mean that the moduli of the arithmetic progression is not equal to 1.} arithmetic progression. The set $A \subseteq \N^{(k)}$ can belong to a non-trivial arithmetic progression if and only if
\begin{displaymath}
\underline{\delta}(A) \leq \max_{p} \frac{\max_a|\{b \in [p] \mid b^k \equiv a \Mod p\}|}{p} = \max_{p} \frac{(k, p-1)}{p} = \max_{p: p-1 | k} \frac{p-1}{p}.
\end{displaymath}

The following theorem shows that the density condition $\underline{\delta}(A) > \max_{p: p-1 | k} \frac{p-1}{p}$ can be obtained if we assume that the number of summands is very large depending on $k$.

\begin{Theorem} \label{theorem_second} Let $k \geq 2$ and $\delta > 0$. Let $A \subseteq \N^{(k)}$ be such that $A$ is not a subset of any non-trivial arithmetic progression and $\underline{\delta}(A) > \delta$. There exists $s = s(k, \delta) \in \N$ such that all sufficiently large natural numbers belongs to the set $sA$.
\end{Theorem}

\subsection{Outline of the proof of Theorem \ref{theorem_main}}

We prove Theorem \ref{theorem_main} using the transference principle, which we introduce in Section \ref{section_transference}. 

Let $f$ be, roughly speaking, a characteristic function of the set $A$ in Theorem $\ref{theorem_main}$. In order for the transference principle to work we need that the function $f$ satisfies three conditions. 1) $f$ needs to satisfy a sufficient mean condition. 2) $f$ has to have a pseudorandom majorant function. 3) $f$ has to satisfy a suitable restriction estimate. We establish these conditions in Sections \ref{section_mean_value_estimate}, \ref{section_pseudorandomness} and \ref{section_restriction_estimate} respectively. 

Both the pseudorandomness condition and the restriction estimate can be dealt with the standard circle method machinery with minor alterations. The mean condition also follows from simple calculations. 

Another main ingredient in the proof of Theorem \ref{theorem_main} is solving the local density version of Waring's problem. Essentially we want to prove that if $A \subseteq Z := \{a \in \Z_{P(w)}^{(k)} \mid (a, P(w)) = 1\}$ and $|A| > \frac{1}{2} |Z|$, then $sA = \Z_{P(w)}$ for some suitably large $s$ depending on $k$, where $P(w) = \prod_{p \leq w} p^k$ and $w \in \N$. We prove this in Section \ref{section_local_version}. This is done using the Chinese remainder theorem, Hensel's lemma and Cauchy-Davenport theorem. 
\\

\noindent \textbf{Remark 1.} The transference lemma (Proposition \ref{proposition_transference}) gives us limition $\delta_A > 2^{-1/k}$. Our result (Theorem \ref{theorem_main}) comes close to this when $k$ is sufficiently large. In particularly for small $k$, we have some density loss because it is not possible to prove the pseudorandomness condition for $f_b$ (for the definition of $f_b$ see (\ref{defn_fb})), when $(W, b) > 1$: If $(W, b) > 1$, we will eventually lose $w$-smoothness of $W$ in calculations, which is crucial for proving the pseudorandomness. There is a way to define $f_b$ so that it satisfies pseudorandomness condition for all $b \in \Z_W^{(k)}$, but this leads to a significantly more difficult local problem, which we were not able to solve.
\\\\
\textbf{Acknowledgments} The author wants to thank Oleksiy Klurman for suggesting this interesting problem. The author thanks Trevor Wooley for showing the alternative way of doing the $\epsilon$-removal that we present in Section \ref{section_restriction_estimate}. The author is also grateful to his supervisor Kaisa Matomäki for many useful discussions. The author thanks the referees for careful reading of the paper and for useful comments. The author is thankful to Victor Wang for pointing out a problem in the proof of Theorem \ref{theorem_second} and suggesting a workaround for the problem. During the work author was supported by Emil Aaltonen foundation.

\section{Notation}

For the rest of the paper we are going to assume that $k \geq 2$ is a fixed integer.

Let $s \in \N$ and $s \geq 2$. For the set $A \subseteq \N$ we define the sumset by
\begin{displaymath}
sA = \{a_1 + \dots + a_s \text{ }|\text{ } a_1, \dots, a_s \in A\}.
\end{displaymath}
For any integers $q, b$, we define the sets
\begin{displaymath}
b + A = \{b\} + A
\end{displaymath}
and
\begin{displaymath}
q \cdot A = \{qa \text{ }|\text{ } a \in A \}.
\end{displaymath}

For finitely supported functions $f,g: \Z \rightarrow \C$, we define convolution $f*g$ by
\begin{displaymath}
f*g(n) = \sum_{a+b=n} f(a)g(b).
\end{displaymath}

For a set $A$, write $1_A(x)$ for its characteristic function. Let $A, B \subseteq [N]$ and $\eta > 0$. We define $S_{\eta}(A, B)$ by
\begin{displaymath}
S_{\eta}(A, B) =  \{n:1_A*1_B(n) \geq \eta N\}.
\end{displaymath}

The Fourier transform of a finitely supported function $f: \Z \rightarrow \C$ is defined by
\begin{displaymath}
\widehat{f}(\alpha) = \sum_{n \in \Z} f(n)e(-n\alpha)
\end{displaymath}
where $e(x) = e^{2\pi i x}$. We will also use notation $e_W(n)$ as an abbreviation for $e(n/W)$.

\indent Let $f: \R \rightarrow \C$ and $g: \R \rightarrow \R_+$. We write
$f = O(g), f \ll g$
if there exists a constant $C > 0$ such that $|f(x)| \leq C g(x)$
for all values of $x$ in the domain of $f$. If $f$ takes only positive values we then define similarly  $f \gg g$ if there exists a constant $C > 0$ such that $f(x) \geq C g(x)$ for all values of $x$ in the domain of $f$. If the implied constant $C$ depends on some contant $\epsilon$ we use notations $O_\epsilon, \ll_\epsilon, \gg_\epsilon$. If $f \ll g$ and $f \gg g$ we write $f \asymp g$. We also write $f = o(g)$ if
\begin{displaymath}
\lim _{x \rightarrow \infty} \frac{f(x)}{g(x)} = 0.
\end{displaymath}
The function $f$ is asymptotic to $g$, denoted $f \sim g$ if 
\begin{displaymath}
\lim _{x \rightarrow \infty} \frac{f(x)}{g(x)} = 1.
\end{displaymath}

We will use notation $\T$ for $\R/\Z$. We also define the $L^p$-norm
\begin{displaymath}
||f||_{p} = \Big(\int_{\T}|f(\alpha)|^pd\alpha \Big)^{1/p}
\end{displaymath}
for function $f: \T \rightarrow \C$.

\section{Transference principle} \label{section_transference}

In this section, we apply the transference principle to prove the transference lemma (Proposition \ref{proposition_transference} below), which we use to prove our main theorem. The idea of the transference principle is to transfer an additive combinatorial result from the integers to a sparse subset of the integers. Particularly these sparse subsets need to be pseudorandom. 

\subsection{The sumset problem in dense settings}

In this subsection, we prove the sumset result, where the sets of the problem are positive density subsets of natural numbers. We later transfer the solution of this dense problem, using the transference principle, to the solution of our sparse problem (the density version of Waring's problem).

We need the following lemma from \cite[Corollary 6.2]{green} that is quantitative version of Cauchy-Davenport theorem.

\begin{Lemma}\label{lemma_davenport} Let $\eta > 0$ and $p$ be a prime. Let $A, B \subseteq \Z_p$ and $|A|, |B| \geq \sqrt{\eta}p$. Then
\begin{displaymath}
|S_\eta(A, B)| \geq \min(p, |A| + |B| - 1) - 3 \sqrt{\eta}p.
\end{displaymath}
\end{Lemma}

Using the previous lemma inductively we prove the following result.

\begin{Lemma}\label{lemma_gen_davenport} Let $p$ be a prime, $s \in \N$, $s \geq 2$, $\epsilon > 2s/p$ and let $B_1, \dots, B_s \subseteq \Z_p$ be such that $\sum_i |B_i| > (1 + \epsilon)p$ and $|B_i| > (\epsilon/s)p$ for all $i \in \{1, \dots, s\}$. Then, for all $n \in \Z_p$, we have 
\begin{displaymath}
1_{B_1}*\dots * 1_{B_s}(n) \gg_{\epsilon, s} p^{s-1}.
\end{displaymath}
\end{Lemma}

\begin{proof}

Let $\eta = \epsilon/6s^2$,
\begin{eqnarray}
R_1 := B_1 \text{ and } R_{i+1} := S_{\eta^2}(R_i, B_{i+1}) \label{R_i}
\end{eqnarray}
for all $i \in \{1, \dots, s-1\}$. Now it follows from Lemma \ref{lemma_davenport} that
\begin{displaymath}
|R_{2}| = |S_{\eta^2}(B_1, B_2)| \geq \min(p, |B_1| + |B_2| - 1) - 3\eta p.
\end{displaymath}
Similarly 
\begin{eqnarray*}
|R_{3}| = |S_{\eta^2}(R_2, B_3)| &\geq & \min(p, |R_2| + |B_3| - 1) - 3\eta p \\
&\geq & \min(p, |B_1| + |B_2| + |B_3| - 2) - 6 \eta p.
\end{eqnarray*}
Repeating this argument inductively, for each $i \in \{1, \dots, s-2\}$, we get that
\begin{displaymath}
|R_{s-1}| \geq \min\Big(p, \sum_{1 \leq i \leq s-1} |B_i| - (s-2)\Big) - 3(s-2)\eta p.
\end{displaymath}
For $n_0 \in \N$ let $N(n_0) := |\{(a, b) \in R_{s-1} \times B_s: a+b \equiv n_0 \Mod p\}|$. We see that
\begin{displaymath}
N(n_0) = |B_s \cap (n_0 - R_{s-1})|  = |B_s \setminus (\Z_p \setminus (n_0 - R_{s-1}))|.
\end{displaymath}
Hence
\begin{eqnarray*}
N(n_0)
&\geq & |R_{s-1}| - ( p - |B_s|)\nonumber\\
&\geq & \min\Big(p + |B_s|, \sum_{i=1}^s |B_i| - (s-2)\Big) - 3(s-2)\eta p- p \nonumber\\
& > & \min\Big(p + (\epsilon/s)p, (1+\epsilon)p - s\Big) - 3s\eta p- p \nonumber\\
& = & \epsilon' p, \label{number_of_ways}
\end{eqnarray*} 
where $\epsilon' = \min(\epsilon/s - 3s\eta,\epsilon - (3s\eta + s/p)) > \epsilon / 4$. Now, for all $n \in \Z_p$, we have that
\begin{eqnarray*}
1_{B_1}* \dots * 1_{B_{s}}(n) 
&\geq & \sums{a + b = n\\a\in R_{s-1}\\b \in B_{s}} 1_{B_1}* \dots * 1_{B_{s-1}}(a)1_{B_{s}}(b)\\
&\geq &  \epsilon'  p \min_{b \in R_{s-1}} 1_{B_1}* \dots * 1_{B_{s-1}}(b)\\
&\geq & \epsilon'  p \min_{b \in R_{s-1}} \sums{i + j = b\\i\in R_{s-2}\\j \in B_{s-1}} 1_{B_1}* \dots * 1_{B_{s-2}}(i)\\
&\geq &  \epsilon'  p \eta^2 p \min_{i \in R_{s-2}}  1_{B_1}* \dots * 1_{B_{s-2}}(i).\\
\end{eqnarray*}
Repeating the last two steps in the previous argument $s-3$ times, it follows that
\begin{displaymath}
1_{B_1}* \dots * 1_{B_{s}}(n) \geq \epsilon'\eta^{2(s-2)}p^{s-1}. \qedhere
\end{displaymath}
\end{proof}

Now we are ready to prove the following sumset lemma. 

\begin{Lemma} \label{lemma_dense_model} Let $\epsilon > 0$, $s \geq 2$ and let $A_1, \dots, A_s \subseteq [N]$ be such that $\sum_i |A_i| > (s(1+ \epsilon)/2)N$ and $|A_i| > (\epsilon / 2) N$ for all $i \in \{1, \dots, s\}$. Then there exists $c(\epsilon, s) > 0$ such that, for all $n \in \Big((1-\frac{\epsilon^2}{16})\frac{sN}{2}, (1 + \frac{\epsilon}{4})\frac{sN}{2}\Big)$, we have
\begin{displaymath}
1_{A_1}*\dots * 1_{A_s}(n) \geq c(\epsilon, s) N^{s-1},
\end{displaymath}
provided that $N$ is sufficiently large depending on $\epsilon$.
\end{Lemma}

\begin{proof} 
Let $p$ be a prime such that $p \in \Big(\frac{(1+\kappa)sN}{2}, \frac{(1+2\kappa)sN}{2}\Big)$, where $\kappa = \epsilon/4$. Such a prime exists by the prime number theorem provided that $N$ is large enough depending on $\epsilon$. For $i \in \{1, \dots, s\}$ define $B_i \subseteq \Z_p$ by $B_i := \{a \Mod p: a \in A_i\}$. We see that 
\begin{displaymath}
\sum_{i=1}^s |B_i| = \sum_{i=1}^s |A_i| > \frac{s (1 + \epsilon)}{2} N = \frac{1 + 4\kappa}{2} sN >  (1+\kappa') p,
\end{displaymath}
where $\kappa' = \frac{2\kappa}{1 + 2\kappa}$. Similarly $|B_i| > (\kappa'/s) p$ for all $i \in \{1, \dots, s\}$. Assuming that $N$ is sufficiently large depending on $\epsilon$, we have that $\kappa' > 2s/p$. Hence it follows from Lemma \ref{lemma_dense_model} that, for any $n \in \Z_p$, 
\begin{displaymath}
1_{B_1}*\dots * 1_{B_s}(n) \gg_{\epsilon, s} p^{s-1} \gg_{\epsilon, s} N^{s-1}.
\end{displaymath}
For each integer $n \in A_1 + \dots + A_s$ we have $n \leq sN < \frac{2}{1 + \kappa} p$. On the other hand, for $n \in \big(\frac{1-\kappa}{1+\kappa}p, p\big)$, we have $p + n > \frac{2}{1 + \kappa} p$. Thus, for $n \in \Big(\frac{1 - \kappa^2}{2}sN, \frac{1 + \kappa}{2}sN \Big)$, we have
\begin{displaymath}
1_{B_1}*\dots * 1_{B_s}(n) = 1_{A_1}*\dots * 1_{A_s}(n)
\end{displaymath}
and the claim follows.
\end{proof}

\subsection{Transference}

In this subsection, we establish the transference lemma, which we will use to prove our main theorem. But first, we introduce some necessary definitions.

\begin{Definition} Let $\eta > 0$ and $N \in \N$. We say that function $f: [N] \rightarrow \R_{\geq 0}$ is $\eta$\textbf{-pseudorandom} if there exists a majorant function $\nu_f$ such that $f \leq \nu_f$ pointwise and $||\widehat{\nu_f} - \widehat{1_{[N]}}||_{\infty} \leq \eta N$.
\end{Definition}

\begin{Definition} Let $q > 1$, $N \in \N$ and $K \geq 1$. We say that function $f: [N] \rightarrow \R_{\geq 0}$ is $q$\textbf{-restricted with constant }$K$ if $||\widehat{f}||_q \leq K N^{1-1/q}$.
\end{Definition}

\begin{Definition} Let $\delta > 0$ and $N \in \N$. We say that function $f: [N] \rightarrow \R$ is $\delta$\textbf{-uniform} if $||\widehat{f}||_\infty \leq \delta N$.
\end{Definition}

Let $N \in \N$, $\delta > 0$ and $f: [N] \rightarrow \R_{\geq 0}$ be a function. Let $T$ be the set of large frequencies of $f$:
\begin{displaymath}
T := \{\gamma \in \T: |\widehat{f}(\gamma)| \geq \delta N\}
\end{displaymath}
We also define a Bohr set using these frequencies:
\begin{displaymath}
B(\delta, N) = \{1 \leq b \leq \delta N: ||b\gamma||_\T < \delta /2 \pi:\forall \gamma \in T\}.
\end{displaymath} 
For the choice of $N, \delta, f$ we define $f_{\delta, N}^*(n) := \E_{a, b \in B}f(n + a - b)$ and $f_{\delta, N}^{unf} := f - f_{\delta, N}^*$. 

Now we can state the following lemma that is the core of the transference principle.

\begin{Lemma} \label{lemma_transfer} Let $\delta > 0$, $N \in \N$ and $K \geq 1$. Let $f: [N] \rightarrow \R_{\geq 0}$ be $\eta$-pseudorandom and $q$-restricted with constant $K$. Then
\begin{enumerate}[label = \rm{(\roman*)}]
\item $0 \leq f_{\delta, N}^*(n) \leq 1 + O_\delta(\eta)$ for all $n \in [N]$
\item $f_{\delta, N}^{unf}$ is $\delta$-uniform
\item $f_{\delta, N}^*$ and $f_{\delta, N}^{unf}$ are $q$-restricted with constant $K$.
\end{enumerate}
\end{Lemma}

\begin{proof}
See the proof of \cite[Lemma 4.3]{matomaki}.
\end{proof}

Next, we prove that the functions $f_1* \dots *f_s$ and $(f_1)_{\delta, N}^* *\dots * (f_s)_{\delta, N}^*$ are in a certain sense close to each other.

\begin{Lemma}\label{lemma_holder_step} Let $\delta > 0$, $\eta > 0$, $N \in \N$ and $K \geq 1$. Let also $s \in \N$, $q \in (s-1, s)$ and, for each $i \in \{1, \dots, s\}$, let $f_i: [N] \rightarrow \R_{\geq 0}$ be a function that is $\eta$-pseudorandom and $q$-restricted with constant $K$. Then, for all $n \in [N]$,
\begin{displaymath}
|f_1* \dots *f_s(n) - (f_1)_{\delta, N}^* *\dots * (f_s)_{\delta, N}^*(n)| \leq  2^s \delta^{s-q} K^q N^{s-1}.
\end{displaymath}

\end{Lemma}

\begin{proof}
Denote $f_i^{unf} = (f_i)_{\delta, N}^{unf}$ and $f_i^* = (f_i)_{\delta, N}^*$ for all $i \in \{1, \dots, s\}$. We see that
\begin{displaymath}
f_1* \dots *f_s(n)  = f_1^* *\dots * f_s^*(n) + \sums{g_i \in \{f_i^*, f_i^{unf}\}\\ \exists i: g_i = f_i^{unf}} g_1* \dots * g_s(n).
\end{displaymath}
Now choose $a = q-s+1 \in (0, 1)$. Let $i \in \{1, \dots, s\}$ be such that $g_i = f_i^{unf}$. Without loss of generality we can assume that $i = 1$. By Hölder's inequality and Lemma \ref{lemma_transfer} we have that

\begin{eqnarray*}
|g_1* \dots * g_s(n)|
&\leq & \int_{\T} |\widehat{g_1}(\gamma) \cdots  \widehat{g_s}(\gamma)| d\gamma \\
&\leq &  ||\widehat{f^{unf}_1}||_{\infty}^{1-a}\int_{\T} |\widehat{f^{unf}_1}|^{a}|\widehat{g_2}(\gamma) \cdots  \widehat{g_s}(\gamma)| d\gamma \\
&\leq &  ||\widehat{f^{unf}_1}||_{\infty}^{1-a}||(\widehat{f^{unf}_1})^a||_{q/a} ||\widehat{g_2}||_q  \cdots ||\widehat{g_s}||_q\\
&=&  ||\widehat{f^{unf}_1}||_{\infty}^{1-a}||\widehat{f^{unf}_1}||_q^{a} ||\widehat{g_2}||_q  \cdots ||\widehat{g_s}||_q\\
&\leq & (\delta N)^{1-a} K^{a} N^{a(1-1/q)} K^{s-1} N^{(s-1)(1-1/q)} \\
& = & \delta^{1-a} K^{s-1+a}  N^{s-1}.
\end{eqnarray*}
Thus 
\begin{displaymath}
\Big| \sums{g_i \in \{f_i^*, f_i^{unf}\}\\ \exists i: g_i = f_i^{unf}} g_1* \dots * g_s(n) \Big| \leq 2^s \delta^{s-q} K^q  N^{s-1}. \qedhere
\end{displaymath}
\end{proof}

Now we are ready to prove the transference lemma which we use to prove our main theorem.

\begin{Proposition} \textbf{(Transference lemma)} \label{proposition_transference}
Let $s \geq 2$, $s - 1 < q < s$, $K \geq 1$ and $\epsilon, \eta \in (0, 1)$. Let $N$ be a natural number and, for each $i \in \{1, \dots, s\}$ let $f_i: [N] \rightarrow \R_{\geq 0}$ be a function that is $\eta$-pseudorandom and $q$-restricted with constant $K$. Assume also that
\begin{equation} 
\mathbb{E}_{n \in [N]} f_1(n) + \dots + f_s(n) > s(1 + \epsilon)/2 \label{eq_*}
\end{equation}
and
\begin{equation}
\mathbb{E}_{n \in [N]} f_i(n) > \epsilon / 2   \label{eq_**}
\end{equation}
for all $i \in \{1, \dots, s\}$. Write $\kappa := \epsilon / 32$. Assume that $\eta$ is sufficiently small depending on $\epsilon, K, q$ and $s$. Then, for all $n \in \Big((1-\kappa^2)\frac{sN}{2}, (1 + \kappa)\frac{sN}{2}\Big)$, we have
\begin{displaymath}
f_1*\dots * f_s(n) \geq c(\epsilon, s) N^{s-1},
\end{displaymath}
where $c(\epsilon, s) > 0$ is a constant depending only on $\epsilon$ and $s$.
\end{Proposition}

\begin{proof} Let $\delta \in (0, \epsilon/8)$ to be chosen later depending on $\epsilon, s, K$ and $q$. Denote $f_i^{unf} = (f_i)_{\delta, N}^{unf}$ and $f_i^* = (f_i)_{\delta, N}^*$ for all $i \in \{1, \dots, s\}$. Write $\lambda := \epsilon/8$ and let $A_i = \{n : f_i^*(n) > \lambda \}$ for all $i \in \{1, \dots, s\}$. By Lemma \ref{lemma_holder_step} we get that
\begin{eqnarray}
f_1*\dots * f_s(n) &\geq & f_1^* *\dots * f_s^*(n) - 2^s \delta^{s-q} K^q N^{s-1} \nonumber \\
&=& \sums{a_1 + \dots + a_s = n\\ a_i \in [N]}f_1^*(a_1) \cdots f_s^*(a_s) - 2^s \delta^{s-q} K^q N^{s-1}\nonumber \\
&\geq & \sums{a_1 + \dots + a_s = n\\ a_i \in A_i}f_1^*(a_1) \cdots f_s^*(a_s) - 2^s \delta^{s-q} K^q N^{s-1}\nonumber \\
&\geq & \lambda ^s \sums{a_1 + \dots + a_s = n}1_{A_1}(a_1) \cdots 1_{A_s}(a_s) - 2^s \delta^{s-q} K^q N^{s-1}\nonumber \\
&\geq & \lambda^s 1_{A_1}*\dots * 1_{A_s}(n) - 2^s \delta^{s-q} K^q N^{s-1} \label{eq_f_convolution_decomposition}.
\end{eqnarray}
For all $i \in \{1, \dots, s\}$, by the definition of $f_i^*$ and Lemma \ref{lemma_transfer} (ii), we get that
\begin{eqnarray*}
\mathbb{E}_{n \in [N]} f_i^*(n) 
&=& \mathbb{E}_{n \in [N]} f_i(n) - \mathbb{E}_{n \in [N]} f_{i}^{unf}(n) \\
&\geq & \mathbb{E}_{n \in [N]} f_i(n) - \delta\\
& > & \epsilon/2 - \delta.
\end{eqnarray*}
By Lemma \ref{lemma_transfer} (i) we see that
\begin{displaymath}
\E_{n \in [N]} f_i^*(n) \leq \frac{1}{N}\sum_{n \in A_i} (1 + O_\delta(\eta)) + \E_{n \in [N]} \lambda.
\end{displaymath}
Thus by (\ref{eq_**})
\begin{displaymath}
(1 + O_\delta(\eta)) |A_i| > (\epsilon/2 - \delta - \lambda) N > (\epsilon/4) N.
\end{displaymath}
Similarly, using (\ref{eq_*}) in place of (\ref{eq_**}), we get that
\begin{displaymath}
s(1+\epsilon)/2 - s\delta \leq \E_{n \in [N]} f_1^*(n) + \dots + f_s^*(n) \leq   \frac{1}{N} \sum_{i=1}^s\sum_{n \in A_i} (1 + O_\delta(\eta)) + s \E_{n \in [N]} \lambda 
\end{displaymath}
and so
\begin{displaymath} 
(1+ O_\delta(\eta))\sum_{i=1}^s |A_i| > (s(1+\epsilon)/2 - s\delta - s\lambda) N > (s(1 + \epsilon/4)/2)N.
\end{displaymath}
We can assume that $\eta$ is small enough in terms of $\epsilon$ and $\delta$, since otherwise the conclusion can be made trivial. Hence $\sum_i |A_i| > (s(1+\lambda)/2)N$ and $|A_i| > (\lambda/2) N$ for all $i \in \{1, \dots, s\}$. 
Let $c'(\lambda, s)$ be the constant in Lemma \ref{lemma_dense_model}. Then the inequality (\ref{eq_f_convolution_decomposition}) and Lemma \ref{lemma_dense_model} imply that
\begin{displaymath}
f_1*\dots * f_s(n) \geq (\lambda^s c'(\lambda, s)- 2^s \delta^{s-q} K^q)N^{s-1}
\end{displaymath}
for all $n \in  \Big((1-\frac{\lambda^2}{16})\frac{sN}{2}, (1 + \frac{\lambda}{4})\frac{sN}{2}\Big)$.
The result now follows by choosing  $\delta$ to be sufficiently small in terms of $\epsilon, s, K, q$.
\end{proof}

In the previous lemma the condition (\ref{eq_*}) is strict:  If $\mathbb{E}_{n \in [N]} f_1(n) + \dots + f_s(n) \leq s/2$, then the sets $A_1, \dots, A_s$ in the proof of Proposition \ref{proposition_transference} can all be subsets of same non-trivial arithmetic progression, which means that also the sumset $A_1 + \dots + A_s$ is subset of a non-trivial arithmetic progression and so $f_1^* * \dots * f_s^*(n) > 0$ is not true for all $n \in \big(\frac{1-\kappa^2}{2}sN, \frac{1+\kappa}{2}sN\big)$.

\section{Proof of the main theorem}

In this section, we will prove Theorem \ref{theorem_main} using the transference lemma (Proposition \ref{proposition_transference}) assuming some lemmas which we will
prove later. We will also prove Theorem \ref{theorem_second}.

\subsection{Definitions}

Let $A \subseteq \N^{(k)}$, $N \in \N$, $w = \log \log \log N$ and
\begin{equation}\label{defn_W}
W := \prod_{p \leq w}p^k.
\end{equation}
Let $b \in [W]$ be such that $b \in \Z_W^{(k)}$. Define $\sigma_W(b) := |\{z \in \Z_{W} \mid z^k \equiv b \Mod W\}|$. Define functions $f_b, \nu_b: [N] \rightarrow \R_{\geq 0}$ by
\begin{equation}  \label{defn_fb}
f_b(n) := \left\{
	\begin{array}{ll}
		\frac{k}{\sigma_W(b)}t^{k-1} & \mbox{if }  Wn + b = t^k \in A \\
		0 & \mbox{otherwise,}
	\end{array}
\right.
\end{equation}
and
\begin{equation} \label{defn_nub}
\nu_b(n) := \left\{
	\begin{array}{ll}
		\frac{k}{\sigma_W(b)}t^{k-1} & \mbox{if } Wn + b = t^k \in \N^{(k)}\\
		0 & \mbox{otherwise}.
	\end{array}
\right.
\end{equation}
Clearly $f_b(n) \leq \nu_b(n)$ for all $n \in [N]$. The purpose of $W$-trick in the definitions of $f_b$ and $\nu_b$ is to make pseudorandomness of $\nu_b$ possible. The normalization of $f_b$ and $\nu_b$ is used to ensure that $\E_{n \in [M]}\nu_b(n) \sim 1$, when $b \in \Z_W^{(k)}$.

Define $Z(q) := \{a \in \Z_q^{(k)} \mid (a, q) = 1\}$. Define also the function $g: [W] \times \N \rightarrow \R_{\geq 0}$ by
\begin{equation} \label{defn_g}
g(b, M) := \E_{n \in [M]}f_b(n).
\end{equation}

For the rest of the paper we will assume the notation of this subsection. 

\subsection{Key lemmas}

We will apply Proposition \ref{proposition_transference} to the function $f_b$. The following three lemmas (to be proven later) show that the function $f_b$ is $\eta$-pseudorandom, $q$-restricted and satisfies the mean condition of Proposition \ref{proposition_transference}. 

\begin{Proposition}\label{proposition_mean_condition} \textbf{(Mean value lemma)} Let $\epsilon \in (0, 1/6)$ and let $N$ be sufficiently large depending on $\epsilon$. Let $\delta_A > (1 - ( 1/2 - 3\epsilon )\mathcal{Z}_k^{-1})^{1/k}$ and $s \geq 16k\omega(k) + 4k + 4$. Then, for all $n \in \Z_W$ with $n \equiv s \Mod{R_k}$, there exist numbers $b_1, \dots, b_s \in Z(W)$ such that $n \equiv b_1 + \dots + b_s \Mod W$, $g(b_i, N) > \epsilon/2$ for all $i \in \{1, \dots, s\}$ and 
\begin{displaymath}
g(b_1, N) + \dots + g(b_s, N) > \frac{s(1+\epsilon)}{2}.
\end{displaymath}
\end{Proposition}

We will prove Proposition \ref{proposition_mean_condition} in Section \ref{section_mean_value_estimate}. The main ingredient in the proof of Proposition \ref{proposition_mean_condition} is a local density version of Waring's problem. We will state and prove this local problem in Section \ref{section_local_version}.

\begin{Proposition} \label{proposition_pseudorandomness_condition} \textbf{(Pseudorandomness)} Let $\alpha \in \T$. Assume that $\sigma_W(b) \not= 0$ and $(b, W) = 1$. Then 
\begin{displaymath}
|\widehat{\nu_b}(\alpha) - \widehat{1_{[N]}}(\alpha)| = o_{k}(N).
\end{displaymath}
\end{Proposition}

We will prove Proposition \ref{proposition_pseudorandomness_condition} in Section \ref{section_pseudorandomness}. The proof uses a standard circle method analysis of major and minor arcs.

\begin{Proposition} \label{proposition_restriction_estimate} \textbf{(Restriction estimate)} Let $s > \max(k^2 + k, 4k)$. Assume that $\sigma_W(b) \not= 0$ and $(b, W) = 1$. Then there exists $q \in (s-1, s)$ such that
\begin{displaymath}
||\widehat{f_b}||_q \ll_k N^{1-1/q}.
\end{displaymath}
\end{Proposition}

We will prove Proposition \ref{proposition_restriction_estimate} in Section \ref{section_restriction_estimate}. The proof is based on Vinogradov's mean value theorem and the $\epsilon$-removal technique.

\subsection{Conclusion}

Now we are ready to prove Theorem \ref{theorem_main} assuming the propositions presented in the previous subsection.

\begin{proof}[Proof of Theorem \ref{theorem_main} assuming Propositions \ref{proposition_mean_condition}, \ref{proposition_pseudorandomness_condition} and \ref{proposition_restriction_estimate}] Let $n_0$ be a large natural number for which $n_0 \equiv s \Mod{R_k}$. Our goal is to prove that $n_0 \in sA$ provided that $n_0$ is sufficiently large. 

Let $N := \lfloor 2n_0 / sW \rfloor $. Choose $\epsilon \in (0, 1/6)$ such that $\delta_A > (1 - ( 1/2 - 3\epsilon )\mathcal{Z}_k^{-1})^{1/k}$. By Proposition \ref{proposition_mean_condition} there exist $b_1, \dots, b_s \in [W]$ such that $n_0 \equiv b_1 + \dots + b_s \Mod W$, $(b_i \mod W) \in Z(W)$, for all $i \in \{1, \dots, s\}$, and the mean conditions (\ref{eq_*}) and (\ref{eq_**}) of Proposition \ref{proposition_transference} hold for the functions $f_{b_1}, \dots, f_{b_s}$. By Propositions \ref{proposition_pseudorandomness_condition} and \ref{proposition_restriction_estimate} also pseudorandomness condition and restriction condition of Proposition \ref{proposition_transference} hold for the functions $f_{b_1}, \dots, f_{b_s}$ for some $q \in (s-1, s)$, $K > 0$ and for any $\eta > 0$. Assume now that $N$ is sufficiently large depending on $\epsilon$ and $\eta$ is sufficiently small depending on $\epsilon, K, q, s$. Then by Proposition \ref{proposition_transference} 
\begin{displaymath}
f_{b_1}* \dots * f_{b_s}(n) > 0,
\end{displaymath}
for all $n \in \big(\frac{1-\kappa^2}{2}sN, \frac{1+\kappa}{2}sN\big)$, where $\kappa = \epsilon/32$.
This means that, for all such $n$,
\begin{equation} \label{eq_Wn+b}
Wn + b_1 + \dots + b_s \in sA.
\end{equation}
Set $n = (n_0 - b_1 \dots - b_s)/W \in \N$. Then $n \sim sN/2$ and so $n \in \big(\frac{1-\kappa^2}{2}sN, \frac{1+\kappa}{2}sN\big)$ provided that $N$ is sufficiently large in terms of $\kappa$. Thus $n_0 \in sA$ by (\ref{eq_Wn+b}).
\end{proof}

We also prove Theorem \ref{theorem_second}. 

\begin{proof}[Proof of Theorem \ref{theorem_second}] For $B \subset \N$ and $s \in \N$ let
\begin{displaymath}
r_B^s(n) := |\{a_1, \dots, a_s \in B \mid n = a_1 + \dots + a_s \}|.
\end{displaymath}
Let $N \in \N$ be sufficiently large and $A' = A \cap [N]$.
By Cauchy-Schwarz inequality and \cite[Theorem 5.7]{nathanson}
\begin{displaymath}
\Big(\sum_{n \in sA'} r_{A'}^s(n)\Big)^2 \leq |sA'|\sum_{n \in sA'} r_{A'}^s(n)^2 \leq |sA'|\sum_{n \in [sN]} r_{\N^{(k)} \cap [N]}^s(n)^2 \ll_{k, s} |sA'| N^{2s/k - 1}
\end{displaymath}
provided that $s > 2^k$. On the other hand 
\begin{displaymath}
\sum_{n \in sA'} r_{A'}^s(n) \geq |A \cap [N]|^s \gg_{\delta, s} N^{s/k}.
\end{displaymath}
Hence
\begin{equation} \label{eq_N_pos_density}
|s (A \cap [N])| > c(k, s, \delta) N
\end{equation}
for all large $N \in \N$ and for some small constant $c(k, s, \delta) > 0$ that depends on $k$, $s$ and $\delta$. 

For $B \subseteq \Z_{\geq 0}$ we define Shnirel'man density
\begin{displaymath}
\sigma(B) := \inf_{N \in \N} \frac{|B \cap [N]|}{N}.
\end{displaymath}

Next, we prove the following claim.

\begin{Claim} Let $B \subseteq \N$ and $\delta' > 0$ such that 
\begin{equation}
\liminf_{N \rightarrow \infty} \frac{|B \cap [N]|}{N} > \delta'. \label{eq_claim_density}
\end{equation}
Assume that $B$ is not a subset of a non-trivial arithmetic progression. Then $uB$ contains two consecutive natural numbers, for some $u \in \N$ large enough depending on $\delta'$.
\end{Claim}

\begin{proof}[Proof of the Claim 1]
 By (\ref{eq_claim_density}) there exists a non-zero $d \in B - B$ with $d = O_{\delta'}(1)$. Let $D = \{n \in B - B \mid n \equiv 0 \Mod d\}$. By (\ref{eq_claim_density}) and the pigeonhole principle there exists an arithmetic progression $P = \{d n + a \mid n \in \N\}$, for some $a \in [d]$, such that 
 \begin{equation} 
\liminf_{N \rightarrow \infty} \frac{|(B \cap P) \cap [N]|}{N} > \delta' / d. \label{eq_claim_density_ap}
\end{equation}
Hence
\begin{equation} 
\liminf_{N \rightarrow \infty} \frac{|D \cap [N]|}{N} > \delta' / d. \label{eq_density_D}
\end{equation}
Let $H = \{n \in \Z_{\geq 0} \mid nd \in D\}$. We see that $0, 1 \in H$. We see by (\ref{eq_density_D}) that $H$ has a positive Shnirel'man density. It now follows by \cite[Theorem 7.7]{nathanson} that there exists $t = O_{\delta'}(1)$ such that 
\begin{equation}
d\Z \subseteq tB - tB \label{eq_AP_in_difference_set}
\end{equation}

Since $B$ is not contained in any non-trivial arithmetic progression, there exist $u \in \N$ with $u = O_d(1)$ and $a_1, \dots, a_u \in \N$ with $a_i + d \in B - B$ for all $i = 1, \dots, u$ such that $(a_1, \dots, a_u, d) = 1$. Hence $a_1, \dots, a_u$ generates $1 \Mod d$ and so for some $v = O_d(1)$ we have 
\begin{displaymath}
1 \in tB - tB + vB - vB.
\end{displaymath}
Thus set $(t+v)B$ contains two consecutive natural numbers. 
\end{proof}

By (\ref{eq_N_pos_density}) we can see that 
\begin{displaymath}
\liminf_{N \rightarrow \infty} \frac{|sA \cap [N]|}{N} > 0.
\end{displaymath}
Since $A$ does not belong to a non-trivial arithmetic progression neither does $sA$. Thus by Claim 1 we obtain that $usA$ contains two consecutive natural numbers for some $u \in \N$ with $u = O_\delta(1)$. Therefore $\sigma(usA - N) > 0$ for some $N \in usA$. Write $A'' = usA - N$. By \cite[Theorem 7.7]{nathanson} there exists $v \in \N$ with $v = O_{\sigma(A'')}(1)$ such that $vA'' = \N$. Therefore all sufficiently large natural numbers belong to the sumset $(vus)A$ and $vus = O_{\delta, k}(1)$.
\end{proof}

\section{Local problem} \label{section_local_version}

In this section, we study the local density version of Waring's problem. This problem is the key new ingredient in solving the density version of Waring's problem. Recall that $Z(q) = \{a \in \Z_q^{(k)} \mid (a, q) = 1\}$. For prime $p$ and $e \geq 1$, we see by \cite[Chapter 4: §2]{ireland}  that
\begin{equation} \label{eq_size_of_Z}
|Z(p^e)| = \frac{\phi(p^e)}{(k, \phi(p^e))}.
\end{equation}
Also, for $n \in \N$, recalling the notation  $\tau(n, p)$ from Section \ref{section_statement_of_results}, we have
\begin{equation} \label{eq_size_of_Z2}
|Z(n)| = \prod_{p | n} \frac{\phi(p^{\tau(n, p)})}{(k, \phi(p^{\tau(n, p)}))}.
\end{equation}
We also note by Fermat's little theorem and the Chinese remainder theorem that if $a \in Z(q)$, then
\begin{equation} \label{wq_mod1}
a \equiv 1 \Mod{(R_k, q)}, 
\end{equation}
where $R_k$ is as in (\ref{defn_R_k}). The congruence in (\ref{wq_mod1}) is the reason why we have the congruence condition in Theorem \ref{theorem_main} as we are restricted to those elements of $A$, which are coprime to $W$.

We will utilize the following definition. 

\begin{Definition} Let $q, s \in \N$. We say that \textit{$(q, s)$} is a \textit{Waring pair} if, for any $A \subseteq Z(q)$ with $|A| > \frac{1}{2}|Z(q)|$, we have $sA = \{a \in  \Z_{q} \mid a \equiv s \Mod{(R_k, q)}\}$.
\end{Definition}

Our aim is to prove the following proposition.

\begin{Proposition} \label{proposition_local_problem} $(W, s)$ is a Waring pair for any $s \geq  8k\omega(k) + 2k + 2$.
\end{Proposition}

We conjecture that $(W, s)$ is a Waring pair for some $s = O(k)$, but we are satisfied with the number of summands being $o(k^2)$, because the restriction estimate ( Proposition \ref{proposition_restriction_estimate}) gives us a lower bound for the number summands that is of order $k^2$.

One of the main reasons why we are able to solve the local problem is the fact that the Waring pairs have multiplicative-like structure. This behaviour is captured in the following lemma. 

\begin{Lemma} \label{lemma_waring_pair_combine} Let $q, r, s, t \in \N$ and $(q, r) = 1$. If $(q, s)$ and $(r, t)$ are Waring pairs, then $(qr, s+t)$ is a Waring pair.
\end{Lemma}

\begin{proof} Let $A \subseteq Z(qr)$ with $|A| > \frac{1}{2} |Z(qr)|$. By the pigeonhole principle there exists a congruence class $a^* \in Z(q)$ such that the set $B := \{b \in A \mid b\equiv a^* \Mod{q}\}$ satisfies $|B| > \frac{1}{2} |Z(r)|$. Let $n \in \Z_{qr}$ be such that $n \equiv s+t \Mod{(R_k, qr)}$. Since $(q, s)$ is a Waring pair, we have that
\begin{align*}
n & \equiv ta^* + a_1 + \dots + a_s \Mod q.
\end{align*}
for some $a_1, \dots, a_s \in A$ (Note that $a \in Z(q)$ implies $a \equiv 1 \Mod{(R_k, q)}$.) Since $(r, t)$ is a Waring pair, we also see that
\begin{align*}
n & \equiv  b_1 + \dots + b_t + a_1 + \dots + a_s \Mod r
\end{align*}
for some $b_1, \dots, b_s \in B$. Hence by the Chinese remainder theorem and definition of $B$
\begin{displaymath}
 n \equiv  b_1 + \dots + b_t + a_1 + \dots + a_s \Mod{qr}.
\end{displaymath}
\end{proof}

We are going to use this lemma to deal separately with $\prod_{\substack{p \leq w \\ p \nmid k}} p^k$ and $\prod_{\substack{ p \leq w \\ p \mid k}} p^k$ parts of $W$.

\subsection{Single moduli}

In this subsection, we study the local problem in $\Z_{p^k}$. For that purpose, we need the following lemma that tells how the elements in $Z(p^k)$ are distributed in certain cosets of $p \cdot \Z_{p^k}$.

\begin{Lemma} \label{lemma_coset_size} Let $p$ be a prime. For all $a \in Z(p)$, we have
\begin{displaymath}
|\{b \in \Z_{p^k}^{(k)} \mid b \equiv a \Mod p\}| = p^{k - 1 - \tau(k, p)}.
\end{displaymath}
\end{Lemma}

\begin{proof}
For $c \in \Z(p)$ set $B(c):= \{b \in \Z_{p^k}^{(k)} \mid b \equiv c \Mod p\}$. For $b, c \in Z(p)$ and $d \in B(c)$, we see that $bd \in B(bc)$. Hence $|B(c)| \leq |B(bc)|$. Since $Z(p)$ is a group, it follows that $|B(b)| = |B(c)|$ for all $b, c \in Z(p)$. Furthermore $|Z(p^k)| = \sum_{b \in Z(p)}|B(b)|$ and so $|B(b)| = |Z(p^k)| / |Z(p)|$ for all $b \in Z(p)$ and the claim follows from (\ref{eq_size_of_Z}).
\end{proof}

We will also need the following generalization of Cauchy-Davenport theorem from \cite[Theorem 1.1]{cochrane}.

\begin{Lemma} \label{lemma_gen_davenport2} Let $n \geq 1$, and $A_1, \dots, A_n$ be finite, nonempty subsets of an abelian group $G$, such that no $A_i$ is contained in a coset of a proper subgroup of $G$. Then
\begin{displaymath}
|A_1 + \dots + A_n| \geq \min\Big(|G|, \Big(\frac{1}{2} + \frac{1}{2n}\Big)\sum_{i=1}^n |A_i|\Big).
\end{displaymath} 
\end{Lemma}
Essentially this means that if $G$ is finite and $A \subseteq G$ satisfies the coset condition, then $A$ is a basis of order $\lceil 2|G|/|A| \rceil - 1$.

Now we can prove the local problem for the prime power moduli.

\begin{Lemma} \label{lemma_local_problem_single} Let $p$ be a prime. Then $(p^k, s)$ is a Waring pair for all $s \geq 8k$. 
\end{Lemma}

\begin{proof} Let $A \subseteq Z(p^k)$ with $|A| > \frac{1}{2} |Z(p^k)|$. If $p-1 | k$, then we see by (\ref{wq_mod1}) that $A \subseteq \{a \in \Z_{p^k} \mid a \equiv 1 \Mod{p^{\eta(p, k)}}\}$, where $\eta(p, k)$ is as in (\ref{defn_eta}). Define 
\begin{displaymath}
A' = \{a \in \Z_{p^{k-\eta(k, p)}} \mid (ap^{\eta(k, p)} + 1 \text{ mod } p^k) \in A\}.
\end{displaymath}
Since $|A'| = |A| > \frac{1}{2} |Z(p^k)| =  \frac{1}{2} p^{k - 1 - \tau(k, p)} \geq \frac{1}{2} p^{k - \eta(k, p)}$ it follows that $A'$ does not belong to any coset of a proper subgroup of $\Z_{p^{k-\eta(k, p)}}$. Hence by Lemma \ref{lemma_gen_davenport2} we get that
\begin{displaymath}
s A' = \Z_{p^{k-\eta(k, p)}}
\end{displaymath}
for all $s \geq 4$.

Similarly if $p - 1 \nmid k$, then $|A| > \frac{1}{2} |Z(p^k)| = \frac{1}{2} p^{k - 1 -\tau(k, p)}\frac{p-1}{(k, p-1)} \geq p^{k - 1 -\tau(k, p)}$. Thus by Lemma \ref{lemma_coset_size} $A$ does not belong to any coset of a proper subgroup of $\Z_{p^{k}}$. Again by Lemma \ref{lemma_gen_davenport2} we get that
\begin{displaymath}
s A = \Z_{p^{k}}
\end{displaymath}
for all $s \geq \lceil 2 |\Z_{p^k}| / |A|\rceil - 1$. By (\ref{eq_size_of_Z}) and the definition of $A$ we see that 
\begin{displaymath}
\lceil  2 |\Z_{p^k}| / |A|\rceil - 1 <  2 \frac{p^k}{\frac{1}{2}\phi(p^k)/(k, \phi(p^k))}  \leq  4k \frac{p}{p-1}  \leq 8k.
\end{displaymath}
\end{proof}

Using Lemmas \ref{lemma_waring_pair_combine} and \ref{lemma_local_problem_single} we can already see that $(W, s)$ is a Waring pair provided that $s \geq \omega(W)8k$, but this is not sufficient as we want to have $s = o(k^2)$. This means that we cannot use Lemma \ref{lemma_waring_pair_combine} too many times.

\subsection{Large moduli}

In this subsection, we deal with the local problem for $k$-coprime part of $W$. First we use Hensel's lemma to reduce the moduli of the problem to be square-free. Then we use a downset idea from \cite[Section 4]{matomaki_2} to simplify the problem.

We start with the moduli reduction argument. 

\begin{Lemma} \label{lemma_power_reduction_trick} Let $e, s \in \N$. Let $q$ be a square-free natural number with $(q, k) = 1$. If $(q, s)$ is a Waring pair, then $(q^e, s + 2)$ is also a Waring pair.
\end{Lemma}

\begin{proof} Let $A \subseteq Z(q^e)$ be any set with $|A| > \frac{1}{2} |Z(q^e)|$ and let $a \in Z(q)$. Then by the Chinese remainder theorem and Hensel's lemma (see e.g. \cite[Proposition 4.2.3]{ireland}) we have that the equation 
\begin{displaymath}
a + b q \equiv x^k \Mod{q^e}
\end{displaymath}
is soluble for all $b \in \Z_{q^{e-1}}$. Hence we can partition $Z(q^e)$ into sets $a + q \Z_{q^{e-1}}$, where $a$ runs through all elements in $\Z(q)$. By the pigeonhole principle we have that for at least one choice of $b \in Z(q)$ it holds that $|H| > \frac{1}{2} q^{e-1}$, where $H = (b + q \Z_{q^{e-1}}) \cap A$. Therefore $2 H = 2b + q  \Z_{q^{e-1}}$. 

Again by the pigeonhole principle there exists an interval $I := (t, (t+1)q]$ for some $t \in [0, q^{e-1}-1]$ such that $|I \cap A| > \frac{1}{2} |Z(q)|$. Since $(q, s)$ is a Waring pair we can now see that
\begin{displaymath}
2 H + s (I \cap A) = \{a \in  \Z_{q^e} \mid a \equiv s + 2 \Mod{(R_k, q)}\}.
\end{displaymath}
\end{proof}

Before we can use the downset idea we need some necessary definitions. Let $n \in \N$ and $a, b \in \Z_n$. We write that $a < b \Mod n$ if and only if there exist $a', b' \in \{0, \dots, n-1\}$ such that $a' < b'$, $a' \equiv a \Mod n$ and $b' \equiv b \Mod n$. Let $q$ be a square-free natural number. For $v \in \Z_q \cong \prod_{p | q} \Z_p$ we define 
\begin{displaymath}
D(v) := \{b \in \Z_q \mid \forall p | q : 0 \leq b \leq v \Mod p\}.
\end{displaymath}
We say that the set $A \subseteq \Z_q$ is a \textit{downset} if $D(v) \subseteq A$ for all $v \in A$. We also say that $u \in \Z_q^*$ is an \textit{upper bound} for the element $a \in \Z_q$ if $a < u \Mod p$ for all $p | q$. We say that $u \in \Z_q^*$ is an \textit{upper bound} for the set $A \subseteq \Z_q$ if $u$ is an upper bound of all elements in $A$.  For $A \subseteq \Z_q$ and $p | q$ define the number of residue classes $\Mod p$ that occur in the set $A$ by $r(A, p) := |\{a \in [p] \mid \exists b \in A : a \equiv b \Mod p \}|$. We define $u(A) \in \Z_q$ such that
\begin{displaymath}
u(A) \equiv r(A, p) \Mod p
\end{displaymath}
for all $p | q$.

The following lemma reveals us how the downsets can be used to analyse the size of sumsets.

\begin{Lemma} \label{lemma_downset} Let $q$ be a square-free natural number. Let $s \in \N$. Let $A_1, \dots, A_s \subseteq \Z_q^*$. Then there exist downsets $A_1', \dots, A_s' \subseteq \Z_q$ such that
\begin{align*}
& |A_i'| = |A_i|, \\
& u(A_i) \text{ is an upper bound for } A_i',
\end{align*}
for all $i \in \{1, \dots, s\}$, and
\begin{displaymath}
|A_1' + \dots + A_s'| \leq |A_1 + \dots + A_s|.
\end{displaymath}
\end{Lemma}

\begin{proof} Let $p | q$ be a prime and write $r = q/p$. For $A \subseteq \Z_q$ and $a \in \Z_r$ define sets $A(a, p), A[a, p], A^{(p)} \subseteq \Z_q$ such that
\begin{align*}
A(a, p)  &:= (\{a\} \times \Z_p) \cap A,\\
A[a, p]  &:= 
\begin{cases}
\{a\} \times \{0, \dots, |A(a, p)| - 1\}, & \text{ if } A(a, p) \not= \emptyset,\\
\emptyset, & \text{ otherwise,}
\end{cases}\\
A^{(p)} &:= \bigcup_{b \in \Z_r} A[b, p].
\end{align*}
In other words the set $A^{(p)}$ has been constructed in a such way that it has the downset property with respect to the coordinate $p$ and it has same number of elements as the set $A$. Clearly $A^{(p)}(a, p) = A[a, p] $. We also define that $\emptyset + A = \emptyset$. We now see that
\begin{align*}
|A_1 + \dots + A_s| & = \sum_{n \in \Z_r} |(A_1 + \dots + A_s)(n, p) |\\
& \geq \sum_{n \in \Z_r} \max_{\substack{a_1, \dots, a_s \in \Z_r\\ \forall i: A_i(a_i, p) \not= \emptyset\\ a_1 + \dots + a_s = n}}\Big|A_1(a_1, p) + \dots + A_s(a_s, p)\Big|.
\end{align*}
Now using the Cauchy-Davenport inequality (\cite[Theorem 5.4]{tao-vu}) we see that
\begin{align*}
|A_1 + \dots + A_s| & \geq \sum_{n \in \Z_r} \max_{\substack{a_1, \dots, a_s \in \Z_r\\ \forall i: A_i(a_i, p) \not= \emptyset\\ a_1 + \dots + a_s = n}} \min\Big(p, |A_1(a_1, p)| + \dots + |A_s(a_s, p)| - (s-1)\Big)\\
& = \sum_{n \in \Z_r} \max_{\substack{a_1, \dots, a_s \in \Z_r\\ a_1 + \dots + a_s = n}} \Big|A_1[a_1, p] + \dots + A_s[a_s, p]\Big|\\
& = \sum_{n \in \Z_r} |(A_1^{(p)} + \dots + A_s^{(p)})(n, p) |\\
& = |A_1^{(p)} + \dots + A_s^{(p)} |.
\end{align*}

Now the sets $A_1^{(p)}, \dots, A_s^{(p)}$ have a downset type property with respect to the $p$-coordinate. Applying the same process to each remaining coordinates $p' | q$ in turn and noticing
that the process does not forget the downsetness of already handled coordinates, we
finally end up with downsets with desired properties.

\end{proof}

Using the previous lemma and simple combinatorial calculations, we can prove the following lemma.

\begin{Lemma} \label{lemma_combinatorial_trick}  Let $q$ be a square-free natural number with $(q, k)=1$. Then $(q, s)$ is a Waring pair for all $s \geq 2k$.
\end{Lemma}

\begin{proof} Let $A \subseteq Z(q)$ with $|A| > \frac{1}{2} |Z(q)|$. For $n \in \N$ set $\sigma(n) := |Z(n)|$. By (\ref{eq_size_of_Z2}) we see that $\sigma$ is a multiplicative function.  Let $u \in \Z_q$ be such that 
\begin{displaymath}
u \equiv \sigma(p) \Mod p
\end{displaymath}
for all $p | q$. By Lemma \ref{lemma_downset} there exists a downset $A' \subseteq Z(q)$ such that $|A| = |A'|$, $u$ is an upper bound for $A'$ and $|sA'| \leq |sA|$ for all $s \geq 1$. Note that $sA'$ is also a downset. 

Now let $S \subseteq \Z_q$ be the set of all elements that have the upper bound $u$. We see that $|S| = \sigma(q)$. We also have $A', u - A' \subseteq S$.  From $2|A'| > |S|$ it follows that
\begin{displaymath}
|\{u = a + b \mid a, b \in A' \}| = |A' \cap (u - A')| = |A' \setminus (S \setminus (u - A'))| \geq |A'| - (|S| - |A'|) > 0.
\end{displaymath}

Hence $u \in 2A'$. Since $2A'$ is a downset, we see that $D(u) \subseteq 2A'$. Because $kD(u) = \Z_q$ we have that $2kA' = \Z_q$.
\end{proof}

From Lemmas \ref{lemma_power_reduction_trick} and \ref{lemma_combinatorial_trick} we now get the following lemma.

\begin{Lemma} \label{lemma_large_moduli} $(\prod_{\substack{p \leq w \\ p \nmid k}} p^k, s)$ is a Waring pair for all $s \geq 2k + 2$.
\end{Lemma}

\subsection{Conclusion}

Combining the results from the previous subsections, we can now solve the local problem.

\begin{proof}[Proof of Proposition \ref{proposition_local_problem}] Using Lemma \ref{lemma_waring_pair_combine} inductively with Lemma \ref{lemma_local_problem_single} to the primes dividing $k$, we get that $(\prod_{\substack{ p \leq w \\ p \mid k}} p^k, 8k\omega(k))$ is a Waring pair. The result now follows from Lemmas \ref{lemma_waring_pair_combine} and \ref{lemma_large_moduli}.
\end{proof}

\section{Mean value estimate} \label{section_mean_value_estimate}

In this section, we will prove the mean condition (Proposition \ref{proposition_mean_condition}) required in the transfence lemma (Proposition \ref{proposition_transference}).

\subsection{Mean value over $g(b, N)$}
In this subsection, we establish a lower bound for $\E_{b\in Z(W)} g(b, N)$, where $g$ is as in (\ref{defn_g}). 

\begin{Lemma}\label{lemma_first_mean} Let $\epsilon \in (0, 1)$. Then
\begin{displaymath}
\E_{b \in \Z_W^{(k)}} g(b, N) \geq (1 - \epsilon) \delta_A^k
\end{displaymath}
provided that $N$ is large enough depending on $\epsilon$.
\end{Lemma}

\begin{proof}

Let $b \in \Z_W^{(k)}$ and write
\begin{displaymath}
\delta_b := \frac{|A \cap (W \cdot [N]+b)|}{|\N^{(k)} \cap (W \cdot[N] + b)|}.
\end{displaymath}
Since $|\N^{(k)} \cap  (W \cdot[N] + b) | \sim \sigma_W(b)(WN)^{1/k}/W$ we have that
\begin{equation} \label{eq_delta_b}
|A \cap (W \cdot [N]+b)| \sim \sigma_W(b)\frac{(WN)^{1/k}}{W} \delta_b.
\end{equation}
Note also that $\sums{t \leq x\\ t \equiv a \Mod n} kt^{k-1} \sim x^k/n$. Hence
\begin{eqnarray}
g(b, N) &=& \frac{1}{N\sigma_W(b)} \sum_{\substack{t^k \leq WN +b\\ t^k \equiv b \Mod W \\ t^k \in A}}kt^{k-1}\nonumber\\
& = & \frac{1}{N\sigma_W(b)} \sums{z \in [W]\\ z^k \equiv b \Mod W} \sum_{\substack{t^k \in A \cap (W \cdot [N] + b) \\ t \equiv z \Mod W }}kt^{k-1} \nonumber\\
& \geq & \frac{1}{N\sigma_W(b)} \sums{z \in [W]\\ z^k \equiv b \Mod W} \sum_{\substack{t \leq W \big\lfloor \frac{|A \cap (W \cdot [N]+b)|}{\sigma_W(b)} \big\rfloor\\ t \equiv z \Mod W}}kt^{k-1} \nonumber\\
&\geq & (1 - o(1)) \frac{1}{WN} \Big(\frac{W}{\sigma_W(b)}\Big)^k |A \cap (W \cdot [N]+b)|^k \nonumber\\
&\geq & (1 - o(1)) \delta_b^k \label{eq_g_b}.
\end{eqnarray}
Since, for any $b \in \Z_W^{(k)}$,
\begin{displaymath}
 \frac{|\N^{(k)} \cap [WN+W-1]|}{|\Z_W^{(k)}|} = |\N^{(k)}  \cap [WN+W-1] \cap (W \cdot \N + b)| + O(1),
\end{displaymath}
we observe that 
\begin{eqnarray}
\delta_A &\leq & \frac{|A \cap [WN+W-1]|}{|\N^{(k)} \cap [WN+W-1]|}\nonumber\\
& = & \E_{b \in \Z_W^{(k)}} \frac{|A \cap [WN+W-1] \cap (W \cdot \N + b)|}{|\N^{(k)}  \cap [WN+W-1] \cap (W \cdot \N + b)| + O(1)}\nonumber\\
& = & (1 + o(1)) \E_{b \in \Z_W^{(k)}} \frac{|A \cap [WN+W-1] \cap (W \cdot \N + b)|}{|\N^{(k)}  \cap [WN+W-1] \cap (W \cdot \N + b)|}\nonumber\\
&=& (1 + o(1))\E_{b \in \Z_W^{(k)}} \delta_b \label{eq_delta_A}.
\end{eqnarray}
Thus by Hölder's inequality, (\ref{eq_g_b}) and (\ref{eq_delta_A})
\begin{displaymath}
\E_{b \in \Z_W^{(k)}} g(b, N) \geq (1 - o(1))\E_{b \in \Z_W^{(k)}} \delta_b^k \geq (1 - \epsilon) \delta_A^k
\end{displaymath}
for any $\epsilon > 0$ provided that $N$ is large enough depending on $\epsilon$. 
\end{proof}

The lower bound in the previous lemma is essentially the best possible as in case $A = \{n^k  \mid n \leq \delta_A (WN)^{1/k}\}$ one has $g(b, \delta_A^k N) \approx 1$, for all $b \in Z(W)$, and so $\E_{b\in Z(W)} g(b, N) \approx \delta_A^k$. 

Using the previous lemma, we can now prove a similar result for $Z(W)$. Recall that 

\begin{equation*} 
\mathcal{Z}_k =  \lim_{m \rightarrow \infty} \frac{|\Z_{P(m)}^{(k)}|}{|\{a \in \Z_{P(m)}^{(k)} \mid (a, P(m)) = 1\}|}.
\end{equation*}

\begin{Lemma} \label{lemma_Z_k_delta} Let $\epsilon > 0$. Let $\mathcal{Z}_k$ be as in (\ref{defn_Zk}). Then
\begin{displaymath}
\E_{b \in Z(W)} g(b, N) > (1 - \epsilon)(\mathcal{Z}_k \delta_A^k- \mathcal{Z}_k + 1)
\end{displaymath}
provided that $N$ is large enough.
\end{Lemma}

\begin{proof} Since $g(b, N) = 1 + o(1)$, we see by Lemma \ref{lemma_first_mean} that
\begin{eqnarray*}
\E_{b \in Z(W)} g(b, N) &=& \frac{|\Z_W^{(k)}|}{|Z(W)|} \E_{b \in \Z_W^{(k)}} g(b, N) - \frac{1}{|Z(W)|}\sums{b \in \Z_W^{(k)} \\ (b, W) > 1} g(b, N)\\
& \geq & (1 - o(1))\mathcal{Z}_k \E_{b \in \Z_W^{(k)}} g(b, N) - (1 + o(1))\mathcal{Z}_k + 1\\
& \geq & (1 - o(1))(\mathcal{Z}_k \delta_A^k- \mathcal{Z}_k + 1). \qedhere
\end{eqnarray*}
\end{proof}

Next we present the following lemma about the size of $\mathcal{Z}_k$.

\begin{Lemma} Let $k > 4$. Then \footnote{Here $\log_2 n = \log n / \log 2$.}
\begin{displaymath}
\frac{\zeta(k)}{\zeta(2k)} \leq \mathcal{Z}_k \leq \frac{\zeta(k - \log_2(2k))}{\zeta(2k - 2\log_2(2k))}.
\end{displaymath}
\end{Lemma}

\begin{proof}  By (\ref{eq_size_of_Z2}) we see that
\begin{displaymath}
\frac{|\Z_W^{(k)}|}{|Z(W)|} = \prod_{p | W} \frac{|Z(p^k)| + 1}{|Z(p^k)|}  = \prod_{p | W} \Big(1 + \frac{1}{|Z(p^k)|}\Big) = \sum_{d | W^{1/k}} \frac{1}{|Z(d^k)|}.
\end{displaymath}
Thus 
\begin{displaymath}
\mathcal{Z}_k = \lim_{w \rightarrow \infty} \frac{|\Z_W^{(k)}|}{|Z(W)|} = \sums{n=1\\ n \text{ square-free}}^\infty \frac{1}{|Z(n^k)|}.
\end{displaymath}
Let $n$ be a square-free natural number. By (\ref{eq_size_of_Z2}) we have $|Z(n^k)| = \prod_{p | n} \frac{p^{k-1}(p-1)}{(k, p^{k-1}(p-1))}$. Since $\omega(n) \leq \log_2 n$ for all $n > 1$, we have that $k^{\omega(n)} \leq n^{\log_2 k}$. Hence $n^{k-1-\log_2 k}\leq |Z(n^k)| \leq n^k$. Now it follows that
\begin{displaymath}
\frac{\zeta(k)}{\zeta(2k)} \leq \mathcal{Z}_k \leq \frac{\zeta(k - \log_2(2k))}{\zeta(2k - 2\log_2(2k))},
\end{displaymath}
because 
\begin{displaymath}
\sums{n=1\\ n \text{ square-free}}^\infty \frac{1}{n^s} = \frac{\zeta(s)}{\zeta(2s)}
\end{displaymath}
for all $s > 1$.
\end{proof}

The previous lemma implies that $\lim_{k \rightarrow \infty} \mathcal{Z}_k = 1$.  Below we have table that illustrates us the convergence of $\mathcal{Z}_k$.

\begin{center}
\begin{tabular}{ c c | c c}
 $k$ & $\mathcal{Z}_k$ &  $k$ & $\mathcal{Z}_k$ \\ \hline 
 2 & 3.279 & 6 & 1.075 \\ 
 3 & 1.493 & 7 & 1.016 \\  
 4 & 1.570 & 8 & 1.062 \\ 
 5 & 1.071 & 9 & 1.004\\    
\end{tabular}
\\
Values of $\mathcal{Z}_k$ for some small values of $k$.
\end{center}

\subsection{Proof of Proposition \ref{proposition_mean_condition}}

For the proof of Proposition \ref{proposition_mean_condition}, we need the following lemma that is essentially a generalized version of the local problem.

\begin{Lemma} \label{lemma_mean_value} Let $f: Z(W) \rightarrow [0, 1)$ be a function satisfying $\E_{b \in Z(W)} f(b) > 1/2$. Let $s \geq 16k\omega(k) + 4k + 4$. Then, for all $n \in \Z_W$ with $n \equiv s \Mod{R_k}$, there exist numbers $b_1, \dots, b_s \in Z(W)$ such that $n \equiv b_1 + \dots + b_s \Mod W$, $f(b_i) > 0$ for all $i \in \{1, \dots, s\}$ and 
\begin{displaymath}
f(b_1) + \dots + f(b_s) > \frac{s}{2}.
\end{displaymath}
\end{Lemma}

\begin{proof} Write $M := \E_{b \in Z(W)} f(b)$. Let $\mu := \max_{b \in Z(W)} f(b)$, $\lambda := 1 - \mu$ and $A := \{b \in Z(W) : f(b) > \lambda\}$. Note that $\mu \geq M > 1/2$ so that $\lambda < 1/2$ and $A$ is non-empty.

We see that
\begin{displaymath}
M \leq \frac{1}{|Z(W)|}\sum_{b \in Z(W) \setminus A}\lambda + \frac{1}{|Z(W)|}\sum_{b \in A} \mu = \lambda\frac{|Z(W)| - |A|}{|Z(W)|} + \mu \frac{|A|}{|Z(W)|}.
\end{displaymath}
Hence
\begin{displaymath}
|A| \geq \frac{M - \lambda}{\mu - \lambda} | Z(W) | = \frac{M + \mu - 1}{2\mu - 1} | Z(W) |  >  M |Z(W)|,
\end{displaymath}
since $M > 1/2$. Thus 
\begin{equation}\label{eq_density_A}
|A | > \frac{1}{2} |Z(W)|.
\end{equation}
Now it follows from Proposition \ref{proposition_local_problem} that 
\begin{displaymath}
s'A = \{a \in \Z_W : a \equiv s' \Mod{R_k}\}
\end{displaymath}
for all $s' \geq 8k\omega(k) + 2k + 2$. 

Now let $b \in Z(W)$ be such that $f(b)= \mu$ and let $s'' \geq s'$. Then, for each $n \in \Z_W$ with $n \equiv s' + s'' \Mod {R_k}$, there exist $b_1, \dots, b_{s'} \in A$ such that
\begin{displaymath}
n - s'' b \equiv b_1 + \dots + b_{s'} \Mod W
\end{displaymath}
and 
\begin{displaymath}
s'' f(b) + f(b_1) + \dots + f(b_{s'}) > s'' \mu + s'\lambda = (s'' - s')\mu + s'(\mu + \lambda) \geq \frac{s'' - s'}{2} + s' = \frac{s' + s''}{2}.
\end{displaymath}
\end{proof}

Using the previous lemma we can now finish the proof of Proposition \ref{proposition_mean_condition}.

\begin{proof}[Proof of Proposition \ref{proposition_mean_condition}] Since the condition $\delta_A > (1 + (3\epsilon - 1/2)\mathcal{Z}_k^{-1})^{1/k}$ is equivalent to $\mathcal{Z}_k \delta_A^k- \mathcal{Z}_k + 1 > 1/2 + 3 \epsilon$, we see by Lemma \ref{lemma_Z_k_delta} that 
\begin{displaymath}
\E_{b \in Z(W)} g(b, N) > (1- \epsilon) (1/2 + 3 \epsilon) > 1/2 + 2\epsilon
\end{displaymath}
provided that $N$ is large enough depending on $\epsilon$. 

For $b \in Z(W)$, define 
\begin{displaymath}
f(b) := \max \Big(0,  \frac{1}{1+\epsilon} (g(b, N) - \epsilon/2)\Big).
\end{displaymath}
Provided that $N$ is large enough depenging on $\epsilon$, we have that $f(b) \in [0, 1)$. We also see that
\begin{displaymath}
\E_{b \in Z(W)} f(b) \geq \frac{1}{1+\epsilon} \E_{b \in Z(W)} (g(b, N) - \epsilon/2) > 1/2.
\end{displaymath}
Hence by Lemma \ref{lemma_mean_value}, for all $n \in \Z_W$ with $n \equiv s \Mod{R_k}$, there exist numbers $b_1, \dots, b_s \in Z(W)$ such that $n \equiv b_1 + \dots + b_s \Mod W$, $f(b_i) > 0$ for all $i \in \{1, \dots, s\}$ and 
\begin{displaymath}
f(b_1) + \dots + f(b_s) > \frac{s}{2}.
\end{displaymath}
By definition of $f$ we have that $g(b_i, N) > \epsilon/2$ for all $i \in \{1, \dots, s\}$ and
\begin{displaymath}
g(b_1, N) + \dots + g(b_s, N) > \frac{s(1+\epsilon)}{2} + \frac{s\epsilon}{2}. \qedhere
\end{displaymath}
\end{proof}

\section{Pseudorandomness condition} \label{section_pseudorandomness}

In this section, we will establish the pseudorandomness of the function $f_b$ (Proposition \ref{proposition_pseudorandomness_condition}). We use the standard circle method machinery to do so.

Let us first introduce the Hardy and Littlewood decomposition. Let 
\begin{equation} \label{defn_Q_T}
Q := N^{\rho} \text{ and } T := N^{1-\rho}
\end{equation}
for $\rho> 0$ to be chosen later. For $a, q \in \N$ and $(a, q) = 1$, write $\mathfrak{M}(q, a) := \{\alpha : |\alpha - \frac{a}{q}| \leq \frac{1}{T} \}$. Let 
\begin{displaymath}
\mathfrak{M} := \bigcup_{\substack{a=0\\(a, q) = 1\\1 \leq q \leq Q}}^{q-1}\mathfrak{M}(q, a).
\end{displaymath}
If $\rho$ is suitably small and $N$ is sufficiently large, then $T > 2Q^2$ and thus the intervals $\mathfrak{M}(q, a)$ are disjoint. Let also $\mathfrak{m} = \mathbb{T} \setminus \mathfrak{M}$. We call $\mathfrak{M}$ major arcs and $\mathfrak{m}$ minor arcs.

From (\ref{defn_nub}) we have that
\begin{eqnarray}
\widehat{\nu_b}(\alpha) &=& \sum_{n}\nu_b(n)e(n\alpha)\nonumber \\
&=& \frac{e_W(-b\alpha)}{\sigma_W(b)} \sum_{\substack{z \in [W]\\ z^k \equiv b \pmod{ W}}} F(\alpha, z), \label{defn_nu_b_composition}
\end{eqnarray}
where
\begin{equation} \label{defn_F}
F(\alpha, z) := \sum_{\substack{t^k \leq WN+b \\ t \equiv z \pmod{W}}}kt^{k-1}e_W(\alpha t^k).
\end{equation}

\subsection{Minor arcs}

In this subsection, we establish Proposition \ref{proposition_pseudorandomness_condition} in the minor arcs using Weyl's inequality.

\begin{Lemma} \label{lemma_minor_arcs} Let $\alpha \in \MinorArcs$.
Then
\begin{displaymath}
\widehat{\nu_b}(\alpha) \ll_{\rho, k} N^{1-\sigma},
\end{displaymath}
for some small $\sigma = \sigma(\rho) > 0$.
\end{Lemma}

\begin{proof}
Let 
\begin{displaymath}
f(X, \alpha, z) = \sum_{\substack{t^k \leq X \\ t \equiv z \pmod{W}}}e_W(\alpha t^k).
\end{displaymath}
Trivially $|f(X, z, \alpha)| \leq X^{1/k}/W$. Let $\lambda \in (0, 1)$ to be chosen later. Using partial summation we get that
\begin{eqnarray*}
F(\alpha, z)
&=& f(WN+b, \alpha, z) k(WN+b)^{1-1/k} - \int_1^{(WN+b)^{1/k}} f(x^k, \alpha, z)k(k-1)x^{k-2}dx \\
&=& f(WN+b, \alpha, z) k(WN+b)^{1-1/k} - \int_{(WN+b)^{(1-\lambda)/k}}^{(WN+b)^{1/k}} f(x^k, \alpha, z)k(k-1)x^{k-2}dx +O((WN)^{1-\lambda}) \\
\end{eqnarray*}

Choose $1 \leq q \leq T$, $(a, q) = 1$ such that $|\alpha - a/q| \leq \frac{1}{qT}$. Since $\alpha \in \mathfrak{m}$, we have $q > Q$. We can write 
\begin{displaymath}
f(X, \alpha, z) = \sum_{\substack{u \leq \frac{X^{1/k} - z}{W}}}e(\alpha W^{k-1} u^k + g(u)),
\end{displaymath}
where $g(u)$ is polynomial with degree at most $k-1$.

Let $q' = \frac{q}{(q, W^{k-1})}$ and $a' = \frac{W^{k-1}a}{(q, W^{k-1})}$. Then $(q', a') = 1$ and
\begin{displaymath}
|\alpha W^{k-1} - a'/q'| \leq \frac{W^{k-1}}{(q, W^{k-1})} \frac{1}{q'T} \leq \frac{W^{k-1}}{(q, W^{k-1})} \frac{1}{q'^2}.
\end{displaymath}
Now by Weyl's inequlity (see e.g. the proof of \cite[Proposition 4.14]{overholt}), for any $\epsilon > 0$,
\begin{displaymath}
f(X, \alpha, z) \ll_{\epsilon, k} \Big(\frac{X^{1/k}}{W}\Big)^{1+\epsilon}\Big(\frac{W^{k-1}}{(q, W^{k-1})}\frac{1}{q'} + \frac{W}{X^{1/k}} + \frac{W^{k-1}}{(q, W^{k-1})} \frac{W^{k-1}}{X^{1-1/k}} + \frac{q' W^k}{X} \Big)^{\sigma},
\end{displaymath}
where $\sigma = \frac{1}{2^{k-1}}$. By (\ref{defn_W}) we have that $W = o(\log N)$. Since $q > Q$, we also see by (\ref{defn_Q_T}) that $q > N^{\rho}$. Thus, for $X \in [(WN)^{1-\lambda}, WN+b]$, we have
\begin{displaymath}
f(X, \alpha, z) \ll_{\epsilon, k} X^{(1+2\epsilon)/k}(N^{-\rho} + X^{-1/k} + X^{1/k-1} + X^{-1} N^{1 - \rho})^{\sigma} \ll  X^{\frac{1 - \sigma' + 2\epsilon}{k}} \ll N^{\frac{(1 - \sigma' + 3\epsilon)}{k}}
\end{displaymath}
for some $\sigma' = \sigma'(\rho) > 0$ provided that $\lambda$ is sufficiently small depending on $\rho$. Hence
\begin{displaymath}
F(\alpha, z)  \ll_{\epsilon, k} N^{1-\sigma''}
\end{displaymath}
for some $\sigma'' = \sigma''(\rho) > 0$ provided that $\epsilon$ is small enough depending on $\sigma'$. The result now follows from (\ref{defn_nu_b_composition}) and (\ref{defn_F}).
\end{proof}

By summing the geometric series (see e.g. \cite[Lemma 4.7]{nathanson}) we see that 
\begin{equation} \label{eq_qeom}
\widehat{1}_{[N]}(\alpha) \ll ||\alpha||^{-1} \ll N^{1-\rho}
\end{equation}
when $\alpha \in \mathfrak{m}$. Hence we get the following lemma.

\begin{Lemma} \label{lemma_minor_arcs_completed} Let $\alpha \in \mathfrak{m}$. Then 
\begin{displaymath}
|\widehat{\nu}_b(\alpha) - \widehat{1}_{[N]}(\alpha)| \ll_\rho N^{1-\epsilon}
\end{displaymath}
for some $\epsilon = \epsilon(\rho) > 0$.
\end{Lemma}

\subsection{Major arcs}

In this subsection, our aim is to prove Proposition \ref{proposition_pseudorandomness_condition} in the major arcs. The result we will prove is the following. 

\begin{Lemma} \label{lemma_major_arcs} Let $\alpha \in \mathfrak{M}$. Assume that $(b, W) = 1$. Then 
\begin{displaymath}
|\widehat{\nu}_b(\alpha) - \widehat{1}_{[N]}(\alpha)| \ll_{k, \epsilon} w^{-1/k+\epsilon} N
\end{displaymath}
for any $\epsilon > 0$ provided that $\rho$ is sufficiently small depending on $k$.
\end{Lemma}

Let us first introduce the following two auxiliary functions that we will use to tackle the pseudorandomness in the major arcs.
\begin{equation} \label{defn_G}
G_b(\alpha, N) := \sum_{\substack{t^k \leq N\\t^k \equiv b \Mod W}}kt^{k-1}e_{W}(\alpha t^k) 
\end{equation}
and 
\begin{displaymath}
V_q(a, b) := \sums{h \Mod{Wq}\\ h^k \equiv b \Mod W} e_{Wq}(ah^k).
\end{displaymath}

The function $G_b(\alpha, N)$ is called the generating function. Our first goal is to prove an approximation lemma for the generating function.

Let
\begin{displaymath}
S(N) := \sum_{\substack{t^k \leq N\\t^k \equiv b \Mod W}}e_{Wq}(a t^k).
\end{displaymath}
We see that
\begin{equation}
S(N) =  \sums{h \Mod{Wq}\\ h^k \equiv b \Mod W} e_{Wq}(ah^k)\sums{t^k \leq N\\ t \equiv h \Mod{Wq}} 1 = V_q(a, b)\frac{N^{1/k}}{Wq} + O(Wq).
\end{equation}

The following lemma approximates the generating function in the rational numbers.

\begin{Lemma}\label{lemma_rational_G} Let $a, q \in \N$. Then
\begin{displaymath}
G_b(a/q, N) = \frac{V_q(a, b)}{Wq}N + O(WqN^{1-1/k})
\end{displaymath}
\end{Lemma}

\begin{proof}
Using partial summation we see that
\begin{eqnarray*}
G_b(a/q, N)
&=& S(N)k N^{1-1/k} - \int_1^N S(t)(k-1)t^{-1/k}dt \\
&=& k \frac{V_q(a, b)}{Wq}N - \frac{V_q(a, b)}{Wq} \int_1^N (k-1)dt + O(WqN^{1-1/k})\\
&=& \frac{V_q(a, b)}{Wq}N + O(WqN^{1-1/k}). \qedhere
\end{eqnarray*}
\end{proof}

Using the previous lemma we can now prove an approximation lemma for the generating function for all real numbers.

\begin{Lemma} \label{lemma_composition_of_G} Let $a, q \in \N$, $\alpha \in \R$ and $\beta = \alpha - a/q$. Then
\begin{displaymath}
G_b(\alpha, N) - \frac{V_q(a, b)}{q}\sums{t \leq N \\ t \equiv b \Mod W} e_{W}(\beta t) = O(Wq N^{1-1/k} + q |\beta| N^{2-1/k}).
\end{displaymath}

\end{Lemma}

\begin{proof}
We can write 
\begin{displaymath}
G_b(\alpha, N) - \frac{V_q(a, b)}{q}\sums{t \leq N \\ t \equiv b \Mod W} e_{W}(\beta t) = \sums{t \leq N\\ t \equiv b \Mod W} u(t) e_W(\beta t),
\end{displaymath}
where
\begin{displaymath}
u(n) = \begin{cases}
kh^{k-1}e_{Wq}(a h^k) -  \frac{V_q(a, b)}{q} & \text{ if } n = h^k\\
-  \frac{V_q(a, b)}{q} 					  & \text{ otherwise.}
\end{cases}
\end{displaymath}
By Lemma \ref{lemma_rational_G} we see that
\begin{eqnarray*}
U(X) &:=& \sums{t \leq X \\ t \equiv b \Mod W} u(t)\\
& = & G_b(a/q, X) - \frac{V_q(a, b)}{q} \sums{t \leq X \\ t \equiv b \Mod W} 1\\
&\ll & WqX^{1-1/k}.
\end{eqnarray*}
Hence, by partial summation,
\begin{eqnarray*}
\sums{t \leq N\\ t \equiv b \Mod W} u(t) e_W(\beta t) &=& e_W(\beta N) U(N) - \int_1^N U(t)\frac{2\pi i \beta}{W} e_W(\beta t)dt\\
&\ll & Wq N^{1-1/k} + q |\beta| N^{2-1/k}. \qedhere
\end{eqnarray*}
\end{proof}

The following lemma tells us that the rational exponential sum $V_q$ vanishes for small values of $q > 1$. This happens because of the $w$-smoothness of $W$. This is also the reason why we use the W-trick in the definition of $f_b$.

\begin{Lemma} \label{lemma_rational_exponential_sum} Let $a, b, q, k \in \N$ be such that $k \geq 2$ and $(a, q) = (b, W) = 1$. Let $\epsilon > 0$. Then
\begin{displaymath}
V_q(a, b) = \begin{cases}
e_W(ab)\sigma_W(b) 						& \text{ if } $q = 1$\\
\sigma_W(b)  O_{\epsilon, k}(q^{1-1/k + \epsilon}) 	& \text{ if } $q > w$\\
0										& \text{ otherwise.}

\end{cases}
\end{displaymath}
\end{Lemma}

\begin{proof}
Follows from \cite[Lemma 21]{salmensuu} and \cite[Theorem]{hua}.
\end{proof}

We record the following consequence of Lemmas \ref{lemma_composition_of_G} and \ref{lemma_rational_exponential_sum} for later use.

\begin{Lemma} \label{lemma_integral_over_G} Let $s > 2k$. Then
\begin{displaymath}
\int_{\mathfrak{M}}|G_b(\alpha, N)|^s d\alpha \ll_{k} \sigma_W(b)^s (N/W)^{s-1}
\end{displaymath}
provided that $\rho$ is small enough depending on $k$.
\end{Lemma}

\begin{proof} Let $\alpha \in \mathfrak{M}(q, a)$ and $\beta = \alpha - a/q$. By the definition of Hardy-Littlewood decomposition in the beginning of Section \ref{section_pseudorandomness} we have that $| \beta | < 1/N^{1-\rho}$ for $\rho > 0$. By \cite[Lemma 4.7]{nathanson} we see that
\begin{displaymath}
\sums{t \leq N\\ t \equiv b \Mod W} e_W(\beta t) \ll \min(N/W, ||\beta||^{-1}).
\end{displaymath}
Thus by Lemmas \ref{lemma_composition_of_G} and \ref{lemma_rational_exponential_sum} we have that 
\begin{displaymath}
G_b(\alpha, N) \ll_{\epsilon, k} \sigma_W(b)q^{\epsilon - 1/k}\min(N/W, ||\beta ||^{-1}) +  N^{1-1/k + 2\rho}
\end{displaymath}
for any $\epsilon > 0$. Hence
\begin{displaymath}
G_b(\alpha, N) \ll_{\epsilon, k}  \sigma_W(b) q^{\epsilon - 1/k}\begin{cases}
||\beta ||^{-1} & \text{ if } \beta \in [W/N, 1/N^{1-\rho}], \\
N/W & \text{ otherwise},
\end{cases}
\end{displaymath}
provided that $\rho$ is small enough depending on $k$. Therefore
\begin{eqnarray*}
\int_{\mathfrak{M}}|G_b(\alpha, N)|^s d\alpha 
& = & \sums{1 \leq q \leq Q\\ 0 \leq a < q\\  (a, q) = 1} \int_{\mathfrak{M}(q, a)} |G_b(\alpha, N)|^s d\alpha \\
& \ll_{\epsilon, k} &  \sigma_W(b)^s\sum_{q \leq Q} q^{ 1 + (\epsilon - 1/k)s} \Big( \int_{0}^{W/N} (N/W)^s d\beta + \int_{W/N}^{1/N^{1-\rho}} ||\beta ||^{-s} d\beta \Big)\\
& \ll & \sigma_W(b)^s (N/W)^{s-1},
\end{eqnarray*}
provided that $\epsilon$ is small enough depending on $k$.
\end{proof}

\subsection{Conclusion}

Now we are ready to finish the proof of Proposition \ref{proposition_pseudorandomness_condition} by tackling the major arc case.

\begin{proof}[Proof of Lemma \ref{lemma_major_arcs}]
By (\ref{defn_F}) and (\ref{defn_G}) we have 
\begin{equation} \label{equation_G_reformed}
G_b(\alpha, WN+b) = \sum_{\substack{z \in [W]\\ z^k \equiv b \pmod{ W}}} F(\alpha, z) 
\end{equation}
Using (\ref{defn_nu_b_composition}) and Lemma \ref{lemma_composition_of_G}  we see that
\begin{equation}\label{eqn_1}
\begin{split}
\widehat{\nu_b}(\alpha)  &= \frac{e_W(-b\alpha)}{\sigma_W(b)} G_b(\alpha, WN +b)\\
 & = \frac{e_W(-b\alpha)}{\sigma_W(b)} \frac{V_q(a, b)}{q}\sums{t \leq WN + b \\ t \equiv b \Mod W} e_{W}(\beta t)+ O(W^2 q N^{1-1/k} + q W^2 |\beta| N^{2-1/k}).
\end{split}
\end{equation}

As $\alpha \in \mathfrak{M}(q, a)$, we have by (\ref{defn_Q_T}) that $q \leq N^\rho$ and $|\beta| \leq N^{\rho - 1}$. Hence the error term in (\ref{eqn_1}) is $O(N^{1-\epsilon'})$ for some $\epsilon' > 0$ provided that $\rho$ is sufficiently small depending on $k$. 

When $q > 1$ it follows from Lemma \ref{lemma_rational_exponential_sum} that
\begin{displaymath}
\widehat{\nu_b}(\alpha) \ll_{\epsilon, k} w^{\epsilon - 1/k} N
\end{displaymath}
for any $\epsilon > 0$. By (\ref{eq_qeom}) we have $\widehat{1_{[N]}}(\alpha) \ll ||\alpha||^{-1} \ll N^{1-\rho}$.

Hence it remains to analyse the case $q=1$ in which case $a=0$ and $\alpha = \beta$. Therefore
\begin{eqnarray*}
\widehat{1_{[N]}}(\alpha) &=& \sum_{n \leq N} e(n\alpha)\nonumber \\
&=& e_W(-b \alpha)\sums{n \leq WN +b \\n  \equiv b \Mod W} e_W(n\alpha),
\end{eqnarray*}
so by (\ref{eqn_1}) and Lemma \ref{lemma_rational_exponential_sum}
\begin{displaymath}
\widehat{\nu_b}(\alpha) = \widehat{1_{[N]}}(\alpha) + O(N^{1-\epsilon'}). \qedhere
\end{displaymath}
\end{proof} 

Proposition \ref{proposition_pseudorandomness_condition} now follows from Lemmas \ref{lemma_minor_arcs_completed} and \ref{lemma_major_arcs} by choosing $\rho$ to be small enough depending on $k$.

\section{Restriction estimate} \label{section_restriction_estimate}

In this section, we will establish the restriction estimate (Proposition \ref{proposition_restriction_estimate}). We do it by using Vinogradov's mean theorem and $\epsilon$-removal technique. The following lemma is a consequence of Vinogradov's mean value theorem. 

\begin{Lemma} \label{lemma_vinogradov_mean} Let $s \geq k(k+1)$ and $\epsilon > 0$.  Then
\begin{displaymath}
||\widehat{f_b}||_s^s \ll_{k, \epsilon} N^{s-1+\epsilon}.
\end{displaymath}
\end{Lemma}

\begin{proof}
Let $X = (WN + b)^{1/k}$ and $t = \frac{k(k+1)}{2}$. We see that
\begin{eqnarray*}
||\widehat{f_b}||_{2t}^{2t} &=& \int_{\T}|\widehat{f_b}(\alpha)|^{2t} d\alpha \\
&=& \int_\T \sum_{\row{n}{2t}}f_b(n_1)\cdots f_b(n_t)\overline{f_b(n_{t+1})}\cdots\overline{f_b(n_{2t})}\\
&& \hspace{1cm}\cdot e(\alpha(n_1 + \dots + n_t - n_{t+1} - \dots - n_{2t})) d \alpha\\
&=& \sum_{\substack{\row{n}{2t}\\ n_1 + \dots + n_t = n_{t+1} + \dots + n_{2t}}}f_b(n_1)\cdots f_b(n_{2t})\\
&\ll_k & X^{2t(k-1)}\sum_{\substack{z_i \leq X \\z_1^k + \dots + z_t^k = z_{t+1}^k + \dots + z_{2t}^k}} 1\\
&=&  X^{2t(k-1)}\int_\T \Big|\sum_{x \leq X }e(\alpha x^k)\Big|^{2t} d\alpha.
\end{eqnarray*}

Let $J_t^{(k)}(X)$ denote the number of integral solutions of the system 
\begin{displaymath}
\nsum{x^i}{t} = x_{s+1}^i + \dots + x_{2t}^i, \text{ } 1\leq i \leq k,
\end{displaymath}
with $1 \leq \row{x}{2t} \leq X$. 

Now by a triangle inequality application (see \cite[Subsection 2.1]{lilian}) and Vinogradov's mean value theorem (\cite[Theorem 1.1]{bourgain}) we have that

\begin{eqnarray*}
\int_\T \Big|\sum_{x \leq X }e(\alpha x^k)\Big|^{2t} d\alpha &\ll_{t,k} & X^{\frac{k(k-1)}{2}}J_{t, k}(X) \\
&\ll_{t, k, \epsilon} & X^{\frac{k(k-1)}{2}} X^{2t - \frac{k(k+1)}{2} + \epsilon} \\
& \ll & X^{2t - k + \epsilon}
\end{eqnarray*}
for all $\epsilon > 0$. Thus 
\begin{displaymath}
||\widehat{f_b}||_{2t}^{2t} \ll_{t, k, \epsilon} X^{2tk - k + \epsilon} \ll  N^{2t-1+\epsilon/k}.
\end{displaymath}
Since $|\widehat{f_b}(\alpha)| \ll N$, it holds, for any $s \geq 2t$, that $||\widehat{f_b}||_s^s \ll_{k, \epsilon} N^{s-1+\epsilon}$.
\end{proof}

Next we introduce the $\epsilon$-removal technique. The $\epsilon$-removal can be done using Bourgain's  strategy from \cite[Section 4]{bourgain_lambda}, but here we use an alternative strategy that the author learned from Trevor Wooley. 

\begin{Lemma} \label{lemma_epsilon_removal} Let $s_0 \geq 1$ be such that
\begin{displaymath}
||\widehat{f_b}||_{s_0}^{s_0} \ll N^{s_0 -1 + \epsilon}
\end{displaymath}
for all $\epsilon > 0$. Then there exists $\gamma \in (0, 1)$ such that, for all $s \geq \max(s_0 + \gamma, 4k + \gamma)$, we have 
\begin{displaymath}
||\widehat{f_b}||_{s}^{s} \ll_k N^{s - 1}.
\end{displaymath}

\end{Lemma}

\begin{proof}
Let $\gamma \in (0, 1)$ to be chosen later and write $s' = s_0 + \gamma$. Let
\begin{equation} \label{eq_open_f}
f'(\alpha) := \frac{1}{\sigma_W(b)}\sums{t^k \in [WN+b] \cap A \\ t^k \equiv b \Mod W} kt^{k-1}e_W(\alpha t^k).
\end{equation}
Since $ |f'(\alpha)| =|\widehat{f_b}(\alpha)|$ it suffices to bound $||f'||_{s'}^{s'}$.

Define $B := \{\alpha \in \T \mid |f'(\alpha)| > N^{1 - 1/s'}\}$ and $I_t : = \int_B |f'(\alpha)|^t d \alpha$, where $t > 0$. Since 
\begin{displaymath}
\int_{\T \setminus B} |f'(\alpha)|^{s'}d\alpha \ll N^{s'-1},
\end{displaymath}
it suffices to show that
\begin{displaymath}
I_{s'} \ll N^{s'-1}.
\end{displaymath}

By (\ref{eq_open_f}) and the Cauchy-Swartz inequality
\begin{eqnarray*}
I_{s'} & = &\frac{1}{\sigma_W(b)}\sums{t^k \in [WN+b] \cap A \\ t^k \equiv b \Mod W} \int_{\alpha \in B} |f'(\alpha)|^{s'-2}f'(-\alpha)k t^{k-1} e_W(\alpha t^k) d\alpha\\
& \leq & \frac{1}{\sigma_W(b)} \Big(\sums{t^k \in [WN+b] \cap A \\ t^k \equiv b \Mod W} kt^{k-1} \Big)^{1/2} \Big(\sums{t^k \in [WN+b] \\ t^k \equiv b \Mod W}\Big|\int_{\alpha \in B}|f'(\alpha)|^{s'-2}f'(-\alpha) k^{1/2}t^{(k-1)/2} e_W(\alpha t^k)\Big|^2\Big)^{1/2}
\end{eqnarray*}
Let 
\begin{eqnarray*}
J & := & \sums{t^k \in [WN+b] \\ t^k \equiv b \Mod W}  \Big| \int_{\alpha \in B}|f'(\alpha)|^{s'-2}f'(-\alpha)k^{1/2}t^{(k-1)/2} e_W(t^k\alpha)\Big|^2\\
& = &\int_{\alpha \in B} \int_{\beta \in B} |f'(\alpha)|^{s'-2} f'(-\alpha) |f'(\beta)|^{s'-2} \overline{f'(-\beta)} \sums{t^k \in [WN+b] \\ t^k \equiv b \Mod W} kt^{k-1} e_W((\alpha - \beta)t^k) d\alpha d \beta \\
& \ll & \int_{\alpha \in B} \int_{\beta \in B} |f'(\alpha)|^{s'-1} |f'(\beta)|^{s'-1} |g(\alpha - \beta)| d\alpha d \beta,
\end{eqnarray*}
where 
\begin{displaymath}
g(\alpha) = \sums{t^k \in [WN+b] \\ t^k \equiv b \Mod W} kt^{k-1} e_W(\alpha t^k).
\end{displaymath}
By 
\begin{displaymath}
\sums{t^k \in [WN+b] \\ t^k \equiv b \Mod W} kt^{k-1} \ll \sigma_W(b)N
\end{displaymath}
we now see that
\begin{equation} \label{eq_Is}
I_{s'} \ll \frac{1}{\sigma_W(b)} (\sigma_W(b) N)^{1/2} J^{1/2}
\end{equation}

 Let $\mathfrak{m}$ and $\mathfrak{M}$ be as in Section \ref{section_pseudorandomness} with $\rho > 0$ to be defined later. We see that
\begin{equation} \label{eq_J_combined}
J \ll J_{\mathfrak{m}} + J_{\mathfrak{M}},
\end{equation}
where, for $M \subseteq \T$, 
\begin{displaymath}
J_{M}  =  \int_{\alpha \in B}  \int_{\substack{\beta \in B\\ \alpha - \beta \in M}} |f'(\alpha)|^{s'-1} |f'(\beta)|^{s'-1} |g(\alpha - \beta)| d\alpha d \beta.
\end{displaymath}

By Lemma \ref{lemma_minor_arcs} we see that, whenever $\alpha \in \mathfrak{m}$,
\begin{displaymath}
g(\alpha) \ll_{\rho, k} \sigma_W(b) N^{1 - \delta},
\end{displaymath}
for some small $\delta = \delta(\rho) > 0$. Since $s_0 > s' - 1$, it follows from definition of $I_t$ that
\begin{displaymath}
I_{s_0} \geq N^{(1-1/s')(s_0 - (s'-1))} I_{s'-1}
\end{displaymath}
and so by assumption $I_{s_0} \ll N^{s_0 - 1 + \epsilon}$, we have
\begin{displaymath}
I_{s'-1} \leq I_{s_0} N^{(1/s' - 1)(s_0 - (s'-1))} \ll_{\epsilon} N^{s'-2 + \epsilon + (s_0 + 1 - s')/s'} \leq  N^{s'-2 + \epsilon + (1 - \gamma)/s'}
\end{displaymath}
for any $\epsilon > 0$. Thus
\begin{equation} \label{eq_J_minor}
J_{\mathfrak{m}} \ll_{\rho, k} \sigma_W(b)I_{s'-1}^2 N^{1 - \delta} \ll_{\epsilon} \sigma_W(b) N^{2s' - 3 + 2\epsilon - \delta + 2(1 - \gamma)/s'} \ll \sigma_W(b) N^{2(s'-1) - 1}
\end{equation}
provided that $2\epsilon - \delta + 2(1-\gamma)/s \leq 0$. This is true if $\epsilon$, $1 - \gamma$ are small enough depending on $\delta$.

Let us now turn to major arcs. Take $\zeta \in (2k, s'/2)$ and choose $h$ such that $\frac{s'}{2 \zeta} + h(1-\frac{1}{\zeta}) = s'-1$. Then by Hölder's inequality
\begin{eqnarray*}
J_{\mathfrak{M}} & \leq & \Big( \int_{\alpha \in B} \int_{\substack{\beta \in B\\ \alpha - \beta \in \mathfrak{M}}}|f'(\alpha)|^{s'} |g(\alpha - \beta)|^\zeta d\alpha d \beta\Big)^{1/(2\zeta)} \\
&& \times  \Big(  \int_{\alpha \in B} \int_{\substack{\beta \in B\\ \alpha - \beta \in \mathfrak{M}}}|f'(\beta)|^{s'}  |g(\alpha - \beta)|^\zeta d\alpha d \beta\Big)^{1/(2\zeta)} \\
&& \times  \Big(  \int_{\alpha \in B} \int_{\substack{\beta \in B}}|f'(\alpha) f'(\beta)|^h d\alpha d \beta\Big)^{1 - 1/\zeta}. \\
\end{eqnarray*}
Note that $h > s'$ since $s' > 2\zeta$. By Lemma \ref{lemma_integral_over_G} and definition of $h$ we have that
\begin{eqnarray}
J_{\mathfrak{M}} \ll_k (I_{s'} \sigma_W(b)^\zeta N^{\zeta-1})^{2/(2\zeta)} I_h^{2(1 - 1/\zeta)} \nonumber
&\ll& (I_{s'} \sigma_W(b)^\zeta N^{\zeta-1})^{1/\zeta} I_{s'}^{2(1 - 1/\zeta)} N^{2(h-s')(1 - 1/\zeta)} \nonumber \\
&=& \sigma_W(b) I_{s'}^{2 - 1/\zeta} N^{(s'-1)/\zeta - 1}, \label{eq_J_major}
\end{eqnarray}
provided that $\rho$ is small enough depending on $k$.
Combining (\ref{eq_Is}), (\ref{eq_J_combined}), (\ref{eq_J_minor}) and (\ref{eq_J_major}) we get that
\begin{displaymath}
I_{s'} \ll_k  \frac{1}{\sigma_W(b)} (\sigma_W(b) N)^{1/2} \Big(\sigma_W(b) N^{2(s'-1) - 1} + \sigma_W(b) I_{s'}^{2 - 1/\zeta} N^{(s'-1)/\zeta - 1}\Big)^{1/2} 
\end{displaymath}
Hence
\begin{displaymath}
I_{s'} \ll_k N^{s'-1}.
\end{displaymath}
\end{proof}

Proposition \ref{proposition_restriction_estimate} now follows from Lemmas \ref{lemma_vinogradov_mean} and \ref{lemma_epsilon_removal}.

\printbibliography

\end{document}